
\input amstex 

\input amssym
\input amssym.def

\magnification 1200
\loadmsbm
\parindent 0 cm


\define\nl{\bigskip\item{}}
\define\snl{\smallskip\item{}}
\define\inspr #1{\parindent=20pt\bigskip\bf\item{#1}}
\define\iinspr #1{\parindent=27pt\bigskip\bf\item{#1}}
\define\einspr{\parindent=0cm\bigskip}

\define\ot{\otimes}

\define\tr{\triangleright}
\define\tl{\triangleleft}

\centerline{\bf Bicrossproducts of algebraic quantum groups}
\bigskip\bigskip
\centerline{\it L.\ Delvaux$^{\text{1}}$, A.\ Van Daele$^{\text{2}}$ and S.H.\ Wang$^{\text{3}}$}
\bigskip\bigskip\bigskip
{\bf Abstract} 
\bigskip 
Let $A$ and $B$ be two algebraic quantum groups (i.e.\ multiplier Hopf algebras with integrals). Assume that $B$ is a right $A$-module algebra and that $A$ is a left $B$-comodule coalgebra. If the action and coaction are matched, it is possible to define a coproduct $\Delta_\#$ on the smash product $A\# B$ making the pair $(A\#B,\Delta_\#)$ into an algebraic quantum group. This result is proven in 'Bicrossproducts of multiplier Hopf algebras' (reference [De-VD-W]) where the precise matching conditions are explained in detail, as well as the construction of this bicrossproduct.
\snl
In this paper, we continue the study of these objects. First, we study the various data of the bicrossproduct $A\# B$, such as the modular automorphisms, the modular elements, ... and obtain formulas in terms of the data of the components $A$ and $B$. Secondly, we look at the dual of $A\# B$ (in the sense of algebraic quantum groups) and we show it is itself a bicrossproduct (of the second type) of the duals $\widehat A$ and $\widehat B$. The result is immediate for finite-dimensional Hopf algebras and therefore it is expected also for algebraic quantum groups. However, it turns out that some aspects involve a careful argument, mainly due to the fact that coproducts and coactions have ranges in the multiplier algebras of the tensor products and not in the tensor product itself.
\snl
Finally, we also treat some examples in this paper. We have included some of the examples that are known for finite-dimensional Hopf algebras and show how they can also be obtained for more general algebraic quantum groups. We also give some examples that are more typical for algebraic quantum groups. In particular, we focus on the extra structure, provided by the integrals and associated objects. In [L-VD] related examples are considered, involving totally disconnected locally compact groups, but as these examples require some analysis, they fall outside the scope of this paper.  It should be mentioned however that with examples of bicrossproducts of algebraic quantum groups, we do get examples that are essentially different from those commonly known in Hopf algebra theory.
\nl
February 2012 ({\it Version} 1.0)
\vskip 1 cm
\hrule
\bigskip\parindent 0.3 cm
\item{$^{\text{1}}$} Department of Mathematics, Hasselt University, Agoralaan, B-3590 Diepenbeek (Belgium). E-mail: Lydia.Delvaux\@uhasselt.be
\item{$^{\text{2}}$} Department of Mathematics, K.U.\ Leuven, Celestijnenlaan 200B (bus 2400),
B-3001 Heverlee (Belgium). E-mail: Alfons.VanDaele\@wis.kuleuven.be
\item{$^{\text{3}}$} Department of Mathematics, Southeast University, Nanjing 210096, China. E-mail: shuanhwang\@seu.edu.cn
\parindent 0 cm
 
\newpage

\bf 0. Introduction \rm
\nl
Let $A$ and $B$ be regular multiplier Hopf algebras. Assume that we have a right action of $A$ on $B$, denoted by $b\tl a$ when $a\in A$ and $b\in B$, making $B$ into a right $A$-module algebra. Also let $\Gamma:A\to M(B\ot A)$ be a left coaction of $B$ on $A$ making $A$ into a left $B$-comodule coalgebra. Finally, assume that the left coaction of $B$ on $A$ and the right action $A$ on $B$ are matched so that the twisted coproduct $\Delta_\#$ on the smash product $A\# B$ is an algebra map. Then the pair $(A\# B,\Delta_\#)$ is again a regular multiplier Hopf algebra.
\snl
The above result has been obtained in our paper entitled {\it Bicrossproducts of multiplier Hopf algebras} ([De-VD-W]) where all the necessary notions together with the matching conditions have been explained in detail. In particular, see Theorem 2.14 and Theorem 3.1 in [De-VD-W].
\nl
If $A$ and $B$ are algebraic quantum groups, i.e.\ when they are regular multiplier Hopf algebras with integrals, it is shown in [De-VD-W] that the bicrossproduct $(A\# B, \Delta_\#)$ is again an algebraic quantum group. A right integral $\psi_\#$ on $A\# B$ is given by 
$$\psi_\#(a\# b)=\psi_A(a)\psi_B(b)$$
when $\psi_A$ and $\psi_B$ are right integrals on $A$ and on $B$ respectively and where $a\in A$ and $b\in B$. See Theorem 3.3 in [De-VD-W].
\snl
In this paper, we continue the previous work with a more detailed study of this last case.
\nl
\it Content of the paper \rm
\nl
In {\it Section} 1 of this paper, we will see what can already be obtained when we only have integrals on $A$. We will show that there exist distinguished multipliers $y,z\in M(B)$ such that $(\iota_B\ot \varphi_A)\Gamma(a)=\varphi_A(a)y$ and $(\iota_B\ot \psi_A)\Gamma(a)=\psi_A(a)z$ for all $a\in A$ where $\varphi_A$ and $\psi_A$ are a left and a right integral on $A$. We will also show that in fact $y=z$. And we will obtain more properties. We start this section however with the dual case where we assume that $B$ has cointegrals. We finish the section with the $^*$-algebra case and an example to illustrate some of the results obtained here.
\snl
In {\it Section} 2,  we will find formulas for the data, canonically associated with an algebraic quantum group (such as the modular element, the modular automorphisms and the scaling constant), for the bicrossproduct $A\# B$ in terms of the data of the components $A$ and $B$.  First we have the formula for the right integral $\psi_\#$. It is simple and was already found in [De-VD-W]. The formula for the left integral $\varphi_\#$ is also given and it involves the multiplier $y$ obtained in the Section 1. From this,  we get a nice and simple expression for the modular element $\delta_\#$. It is proven that $\delta_\#=\delta_A\# y^{-1}\delta_B$ where $\delta_A$ and $\delta_B$ are the modular elements for $A$ and $B$ respectively and where $y$ is the multiplier in $M(B)$, defined and characterized in Section 1. For the scaling constant $\tau_\#$ we simply have $\tau_\#=\tau_A\tau_B$, where $\tau_A$ and $\tau_B$ are the scaling constants for $A$ and $B$ respectively. The formulas for the modular automorphisms $\sigma_\#$ and $\sigma'_\#$ are also found, but are not at all trivial. The main reason for this is that in general, we can not expect that the right integral on $B$ is invariant for the action of $A$. In some cases it will be relatively invariant (e.g.\ when $B$ is a right $A$-module {\it bi}-algebra). But in general this need not be the case. 
\snl
As in the first section, also here we will see what happens in the case of $^*$-algebras. We also look at the example given at the end of the previous section. For this example we find the modular automorphisms explicitly and see that they are non-trivial in many cases.
\snl
In {\it Section} 3 we obtain the main result about duality. We show that the right action of $A$ on $B$ can be dualized to a right coaction of the dual $\widehat A$ of $A$ on the dual $\widehat B$ of $B$ in a natural way. Also the left coaction of $B$ on $A$ is dualized to a left action of $\widehat B$ on $\widehat A$. This action and  coaction make $(\widehat A,\widehat B)$ into a matched pair (of the second type as in Theorem 2.16 of [De-VD-W]). Moreover, it is shown that the tensor product pairing between the tensor product $A\ot B$ on the one hand and the tensor product $\widehat A\ot \widehat B$ on the other hand, induces an isomorphism of the dual of the smash product $A\# B$ (in the sense of algebraic quantum groups) with the smash product $\widehat A\# \widehat B$. This is the main result of this section, found in Theorem 3.7. The difficulties we encounter in this section are not present in the case of the duality for finite-dimensional Hopf algebras. In Section 3, we focus on these problems, rather than on the explicit (and expected) formulas.
\snl
At the end of this section, we also  briefly treat the $^*$-case.
\snl
In {\it Section} 4 we focus on examples and special cases. We begin with the motivating example, coming from a matched pair of (not necessarily finite) groups as it can be found already in [VD-W]. See also the examples at the end of Section 2 and Section 3 in [De-VD-W]. We are very brief here as this example is only included for completeness. We focus on the integrals and the modular data.
\snl
Next we consider again the example that we considered already at the end of Section 1 and Section 2 of this paper. Here we complete the treatment of this example. In particular, we see how the results of Section 3 apply to this case. We find that, just as in the case of Hopf algebras, the bicrossproduct is isomorphic with the tensor product of the two components. This allows us to verify the (non-trivial) formulas for various data, as well as formulas related with duality.
\snl
We treat the special case where the right action of $A$ on $B$ is trivial in this section, as well as the dual case where the left coaction of $B$ on $A$ is trivial. In the first case, we find that the coaction $\Gamma$ must be an algebra map (so that $A$ is a left $B$-comodule {\it bi}-algebra). In the dual case, we find that $B$ is forced to be a right $A$-module {\it bi}-algebra. There are some more consequences however and later in this section, we also consider generalizations of these special cases where we do no longer assume that the action, resp.\ the coaction is trivial, but where we have a comodule {\it bi}-algebra, resp.\ a module {\it bi}-algebra. In particular, we look again at the case of a matched pair of groups and find such an example in this setting already.
\snl
Finally, in {\it Section} 5 we draw some conclusions and discuss very briefly further research.
\nl
\it Notations, conventions and basic references \rm
\nl
An important notion in the theory is that of the multiplier algebra of an algebra (with a non-degenerate product). We refer to the basic works on multiplier Hopf algebras (see further). And an important tool for working with multiplier Hopf algebras is the possibility of extending linear maps from the algebra to the multiplier algebra. This happens in various situations. If e.g.\ we have a homomorphism $\gamma:A\to M(B)$ for algebras $A$ and $B$, it has a unique  extension to a unital homomorphism, usually still denoted by $\gamma$ from $M(A)$ to $M(B)$, provided the original homomorphism is non-degenerate (that is when $\gamma(A)B=B\gamma(A)=B$). This property can e.g.\ be used to formulate coassociativity of the coproduct. It can also be used to extend the counit and the antipode to the multiplier algebra. In the case of the antipode, one has to consider the algebra with the opposite product. In a similar way, also actions, pairings, etc.\ can be extended, up to a certain level, to the multiplier algebras. We will mention this clearly when we use these extensions for the first time in this paper, but we refer to the literature (see further) for details.
\snl
We will freely use Sweedler's notation for multiplier Hopf algebras, not only for coproducts, but also for coactions. But {\it we will  not do this for the coproduct on the bicrossproduct}, except in a few cases where we will mention this very clearly. This would become very confusing as e.g.\ the coproduct on $A$, as a subalgebra of the multiplier algebra $M(A\# B) $ of the smash product $A\# B$ is not the same as the original copoduct.
\snl
When $A$ is a multiplier Hopf algebra, we use $\Delta_A$, $\varepsilon_A$, and $S_A$ to denote the comultiplication, the counit and the antipode of $A$. Occasionally, and when no confusion is possible, we will drop the index, but most of the time however, we will write the index for clarity. 
\snl
When $A$ and $B$ are regular multiplier Hopf algebras and when $\tl$ is a right action of $A$ on $B$, we will in general write $AB$ for the smash product $A\# B$ and consider it as the algebra 'generated' by $A$ and $B$ subject to the commutation rules $ba=\sum_{(a)}a_{(1)}(b\tl a_{(2)})$. With this convention, the coproduct $\Delta_\#$ on $AB$ is given by the formula 
$$\Delta_\#(ab)=\sum_{(a)}(a_{(1)}\ot 1)\Gamma(a_{(2)})\Delta_B(b)$$ 
when $a\in A$ and $b\in B$.
\snl
When $A$ and $B$ are algebraic quantum groups (as will be the case in most of this paper), we will systematically use $C$ for the dual $\widehat A$ and $D$ for the dual $\widehat B$. This makes the notations consistent with the ones used in [De-VD-W]. 
\snl
When we have a right action of $A$ on $B$ and a left coaction of $B$ on $A$ that are compatible (in the sense that the conditions of Theorem 2.14 of [De-VD-W] are fulfilled), we call $(A,B)$ a matched pair (of the first type). We will show that for algebraic quantum groups, the adjoints of the action and the coaction, mentioned earlier, make $(C,D)$ into a matched pair of the second type. By this we mean that the dual conditions, as formulated in Theorem 2.16 of [De-VD-W] are fulfilled.
\snl
In order to avoid the use of too many different notations, we will use mostly the same notations for concepts associated with the pair $(A,B)$ and with the pair $(C,D)$. We refer to Remark 3.4 in Section 3. On the other hand, we will be very systematic in the use of letters $a,b,c,d$ for elements in resp.\ $A,B,C,D$ so that this practice will not cause too much confusion.
\snl
There appears to be some confusion about the use of the terms {\it bicrossproduct} and {\it double cross product}. We use the term bicrossproduct for the construction with a right action of $A$ on $B$ and a left coaction of $B$ on $A$ (in the sense of Majid, see Section 6.2 of [M]). In this case, both the product and the coproduct is twisted. We would use the term double cross product when we have a right action of $A$ on $B$ and a compatible left action of $B$ on $A$ in which case the product is (double) twisted, but not the coproduct. See Section 7.2 in [M].
\snl 
We refer to [VD1] for the theory of multiplier Hopf algebras and to [VD2] for the theory of algebraic quantum groups, see also [VD-Z]. For the use of the Sweedler notation in this setting, we refer to [VD4]. Finally, for pairings of multiplier Hopf algebras, the main reference is [Dr-VD]. In [De-VD2] relations between data for $(A,\Delta)$ and its dual $(\widehat A,\widehat\Delta)$ are explicitly given. An important reference for this paper is of course [De-VD-W] where the theory of bicrossproducts of multiplier Hopf algebras is developed. Finally, for some related constructions with multiplier Hopf algebras, see [De-VD1] and [Dr-VD].
\nl
\bf Acknowledgements \rm
\nl
We would like to thank an anonymous referee for his/her comments on an earlier version of this material. This led us to rewrite the paper drastically, resulting in two parts with [De-VD-W] as the first part dealing with bicrossproducts and this part where also integrals are considered. 
\nl\nl

\bf 1. Actions and coactions, integrals and cointegrals \rm
\nl
Consider two regular multiplier Hopf algebras $A$ and $B$. Let $\tl$ be a right action of $A$ on $B$  making $B$ into a right $A$-module algebra. Let $\Gamma$ be a left coaction of $B$ on $A$ making $A$ a left $B$-comodule coalgebra. Assume that the action and coaction make of $(A,B)$ a matched pair as discussed already in the introduction.
\snl
In this section, {\it we will consider integrals and cointegrals} for the components $A$ and $B$ and we will investigate the relation with the action and the coaction.
\nl
{\it The case where $B$ has cointegrals}
\nl
We will not be able to deduce very much however, but fortunately,  this case is not so important. We mainly include it for motivational purposes. In particular, we will use this case to explain what kind of problems we encounter.
\snl
Recall that an element $h\in B$ is called a {\it left cointegral} if it is non-zero and if $bh=\varepsilon(b)h$ for all $b\in B$. Similarly, a non-zero element $k\in B$ satisfying $kb=\varepsilon(b)k$ for all $b$ is called a right cointegral. Cointegrals are unique, up to a scalar, if they exist. And of course, if $h$ is a left cointegral, then $S(h)$ will be a right cointegral. 
\snl
It is not so difficult to obtain the following result.

\inspr{1.1} Proposition \rm
Assume that $h$ is a left cointegral and that $k$ is a right cointegral in $B$. Then there exist homomorphisms $
\rho, \eta:A\to \Bbb C$ such that 
$$h\tl a=
\rho(a)h \qquad\qquad \text{and} \qquad\qquad k\tl a = \eta(a)k$$
for all $a\in A$.
\snl \bf Proof\rm:
Let $a\in A$ and $b\in B$. We have
$$\align b(h\tl a) &=\sum_{(a)}((b\tl S(a_{(1)}))h)\tl a_{(2)}\\
                   &=\sum_{(a)}\varepsilon_B(b\tl S(a_{(1)}))(h\tl a_{(2)}).
\endalign$$
We know that $\varepsilon_B(b\tl a')=\varepsilon_B(b)\varepsilon_A(a')$ for all $a'\in A$ and $b\in B$. See a remark following the proof of Theorem 2.14 in [De-VD-W]. So
$$\align b(h\tl a) &=\sum_{(a)}\varepsilon_B(b)\varepsilon_A(S(a_{(1)}))(h\tl a_{(2)})\\
                   &=\varepsilon_B(b)(h\tl a).
\endalign$$
It follows that $h\tl a$ is again a left cointegral (if it is non-zero) and so it is a scalar multiple of $h$. Therefore there is a linear map $\rho:A \to \Bbb C$ defined by $h\tl a=\rho(a)h$. 
\snl
Because $h\tl(aa')=(h\tl a)\tl a'$, we also have $\rho(aa')=\rho(a)\rho(a')$ for all $a,a'\in A$.
\snl
A similar argument will work for the right cointegral and we get also a homomorphism $\eta:A\to \Bbb C$, defined by $k\tl a=\eta(a)k$ for all $a\in A$.
\hfill $\square$
\einspr

In the equations above, all expressions can be covered because the action of $A$ on $B$ is assumed to be unital, see e.g.\ [Dr-VD-Z].  
\snl
We want to make a few remarks before we continue.

\inspr{1.2} Remark \rm
i) We have used that $\varepsilon_B(b\tl a)=\varepsilon_B(b)\varepsilon_A(a)$ for all $a\in A$ and $b\in B$. This property is obvious when $B$ is a right $A$-module {\it bi}-algebra, i.e.\ when also $\Delta_B(b\tl a)=\sum_{(a)(b)}(b_{(1)}\tl a_{(1)})\ot (b_{(2)}\tl a_{(2)})$.  This special case is treated e.g.\ in [De1] where the result is obtained in Lemma 1.5. The covering for the right hand side of this equation is illustrated in Definition 1.4 in [De1].
\snl
ii) If $(A,B)$ is a matched pair as explained before, then $B$ will be a right $A$-module {\it bi}-algebra if the coaction $ \Gamma$ of $B$ on $A$ is trivial (see Proposition 4.13 in Section 4). If this is not the case, it is still true that $\varepsilon_B(b\tl a)=\varepsilon_B(b)\varepsilon_A(a)$ for all $a\in A$ and $b\in B$, but it is more complicated to obtain this result. In any case, it will not be sufficient to require only that $B$ is a right $A$-module algebra. There is needed some relation between the action of $A$ on $B$ and the coproduct $\Delta_B$ on $B$ because we want to say something about the behavior of the counit $\varepsilon_B$ w.r.t.\ the action. 
\hfill $\square$
\einspr

Furthermore, {\it still in the case of an $A$-module bi-algebra $B$}, it is also not so hard to show the following.

\inspr{1.3} Proposition \rm
If $A$ and $B$ are as before and if moreover $B$ is a right $A$-module bi-algebra, then $\rho=\eta$ where the homomorphisms $
\rho$ and $\eta$ are as defined in Proposition 1.1.

\snl\bf Proof\rm: 
We have that in this case, $S_B(b\tl a)=S_B(b)\tl a$ for all $a\in A$ and $b\in B$ (see again Lemma 1.5 in [De1]). Then, with $b=h$ where $h$ is a left cointegral in $B$, we obtain $\rho(a)S(h)=S(h)\tl a=\eta(a)S(h)$ because $S(h)$ is a right cointegral. It follows that $\rho(a)=\eta(a)$ for all $a$.
\hfill $\square$\einspr 

In the {\it general case}, when $B$ is not assumed to be an $A$-module bi-algebra, it is still true that $\rho=\eta$. The proof however is far more complicated because there seems to be no simple formula for $S_B(b\tl a)$ in this case. We will obtain the result later by duality and we will come back to this problem then (see Corollary 3.8 in Section 3). 
\nl
{\it The case where $A$ has integrals} 
\nl 
In this case, we have the following dual version of Proposition 1.1.

\inspr{1.4} Proposition \rm 
Assume that $\varphi$ is a left integral on $A$. Then there is a unique element $y\in M(B)$ such that 
$$(\iota\ot\varphi)\Gamma(a)=\varphi(a)y$$
for all $a\in A$.

\snl \bf Proof\rm:
First recall that $\Gamma(a)(b\ot 1)$ and $(b\ot 1)\Gamma(a)$ are elements in $B\ot A$ for all $a\in A$ and $b\in B$. Therefore, $(\iota\ot\omega)\Gamma(a)$ is well-defined in $M(B)$ for all $a\in A$ and any linear functional $\omega$ on $A$. Also recall that $\Delta_\#(a)=\sum_{(a)}(a_{(1)}\ot 1)\Gamma(a_{(2})$. Therefore, we can define a multiplier $(\iota\ot\omega)\Delta_\#(a)\in M(AB)$ for all $a\in A$ and a linear functional $\omega$ on $A$. 
\snl
Now, take $a\in A$ and consider the equality
$$(\iota\ot\Delta_A)\Gamma(a)=\sum_{(a)}\Gamma_{12}(a_{(1)})\Gamma_{13}(a_{(2)}).$$
Recall that this equality is one of the conditions and has been discussed in Section 1 (see Definition 1.9 and the Remark 1.10) of [De-VD-W]. If we apply this equation to the second leg of $\Delta_A(a)$, we find the formula
$$(\iota\ot\Delta_A)\Delta_\#(a)=\sum_{(a)}(\Delta_\#(a_{(1)})\ot 1)\Gamma_{13}(a_{(2)}).$$
We can cover this formula if we multiply with an element in $AB\ot AB$ (or even an element in $AB\ot A$) in the first two factors. Then we can apply $\varphi$ on the third factor of this equality. We get, using the left invariance of $\varphi$, that
$$(\iota\ot\varphi)\Delta_\#(a)\ot 1 
     =\sum_{(a)}\Delta_\#(a_{(1)})((\iota\ot\varphi)\Gamma(a_{(2)})\ot 1). 
$$
Now, we multiply this equation with $\Delta_\#(p)$ from the left and $q\ot 1$ from the right where $p,q\in AB$. This yields the equality
$$\Delta_\#(p)((\iota\ot\varphi)(\Delta_\#(a)(q\ot 1))\ot 1) 
= \sum_{(a)}\Delta_\#(pa_{(1)})((\iota\ot\varphi)(
\Gamma(a_{(2)})(q\ot 1))\ot 1).$$
From the injectivity of the map $p'\ot q'\mapsto \Delta_\#(p')(q'\ot 1)$ on $AB\ot AB$, it then follows that
$$p \ot (\iota\ot\varphi)(\Delta_\#(a)(q\ot 1)) 
= \sum_{(a)}pa_{(1)} \ot ((\iota\ot\varphi)(
\Gamma(a_{(2)})(q\ot 1)))$$
in $AB\ot AB$. We cancel $p$ and $q$ and rewrite this equality as
$$\sum_{(a)} 1\ot a_{(1)}\ot (\iota\ot\varphi)
\Gamma(a_{(2)}) 
= \sum_{(a)}a_{(1)} \ot 1 \ot(\iota\ot\varphi)
\Gamma(a_{(2)})$$
in $M(A\ot A\ot B)$. Then apply $\varepsilon_A$ to the first leg. We get
$$\sum_{(a)} a_{(1)}\ot (\iota\ot\varphi)
\Gamma(a_{(2)}) 
= 1 \ot(\iota\ot\varphi)
\Gamma(a).$$
It follows that $a\mapsto (\omega(b\,\cdot\,)\ot\varphi)\Gamma(a)$ is a left invariant functional on $A$ for all linear functionals $\omega$ on $B$ and all elements $b\in B$. By uniqueness of left invariant functionals, we have a number $c$, dependent on $\omega$ and on $b$, so that 
$(\omega(b\,\cdot\,)\ot\varphi)\Gamma(a)=c\varphi(a)$ for all $a$. It follows easily that the multiplier $(\iota\ot\varphi)
\Gamma(a)$ has the form $\varphi(a)y$ with $y\in M(B)$. This completes the proof of the proposition.
\hfill $\square$ \einspr

Observe that covering the equations in the proof is non-trivial but no problem as can be seen above, where we multiplied with elements $\Delta_\#(p)$ and $q\ot 1$ for $p,q\in AB$. Similar coverings have to be used at various other places in the proof.
\snl
Remark also that the proof becomes much more easy when the algebras $A$ and $B$ have an identity (i.e.\ when they are Hopf algebras). Indeed, then we can simply use the injectivity of the map $T$, defined on $A\ot B$ by 
$T(a\ot b)=
\Gamma(a)(b\ot 1)$, in an early stage of the proof to get the equation 
$$\sum_{(a)} a_{(1)}\ot (\iota\ot\varphi)
\Gamma(a_{(2)}) 
= 1 \ot(\iota\ot\varphi)
\Gamma(a).$$
\snl
A similar argument would {\it not work} with the right integral. However, the next proposition shows that a similar formula holds for the right integral, actually with the same multiplier $y$, just as in the dual case, proven in Proposition 1.3 only for a right module {\it bi}-algebra.
  
\inspr{1.5} Proposition \rm
Assume that $\psi$ is a right integral on $A$. Then also $(\iota \ot \psi)\Gamma(a)=\psi(a)y$ for all $a$ in $A$ where $y$ is the multiplier in $M(B)$ obtained in Proposition 1.4.

\snl\bf Proof\rm:
We start as in the proof of the previous proposition with an element $a\in A$ and the equality
$$(\iota\ot\Delta_A)\Gamma(a)=\sum_{(a)}\Gamma_{12}(a_{(1)})\Gamma_{13}(a_{(2)}).$$
Now, we apply the antipode $S_A$ on the third factor and multiply to obtain
$$ \sum_{(a)}\Gamma(a_{(1)})(\iota\ot S_A)\Gamma(a_{(2)})
                 =(\iota_B\ot \varepsilon_A)
\Gamma(a) \ot 1_A = \varepsilon_A(a) 1_B\ot 1_A.$$
Replace in this formula $a$ by $a_{(2)}$ and then multiply with $a_{(1)}$ from the left in the first factor to get
$$ \sum_{(a)}\Delta_\#(a_{(1)}) (\iota_B\ot S_A)\Gamma(a_{(2)}) 
   =\sum_{(a)}\varepsilon(a_{(2)})a_{(1)}\ot 1_A =a\ot 1_A.$$
Again replace $a$ by $a_{(2)}$ and now multiply with $\Delta_\#(S_A(a_{(1)}))$. Because we have
$$\align \sum_{(a)}\Delta_\#(S_A(a_{(1)})a_{(2)})(\iota\ot S_A)\Gamma(a_{(3)})
&=\sum_{(a)}\Delta_\#(\varepsilon_A(a_{(1)})1_A)(\iota\ot S_A)\Gamma(a_{(2)})\\
&=(\iota\ot S_A)\Gamma(a)\endalign$$
we get 
$$(\iota\ot S_A)\Gamma(a)=\sum_{(a)}\Delta_\#(S_A(a_{(1)}))(a_{(2)}\ot 1_A).$$
Finally, we apply $\varphi$ to the second leg of this formula. When we assume that $\varphi\circ S_A=\psi$, we find
$$\align (\iota\ot\psi)\Gamma(a)
    &=\sum_{(a)}S_A(a_{(2)})(\iota\ot\varphi)\Gamma(S_A(a_{(1)})) a_{(3)}\\
    &=\sum_{(a)}S_A(a_{(2)})\varphi(S_A(a_{(1)}))ya_{(3)}\\
\endalign$$ 
and because $\varphi\circ S_A$ is right invariant, we find that this last expression is equal to $\varphi(S_A(a))y=\psi(a)y$.
\hfill $\square$\einspr

Remark that we obtain a formula for $(\iota\ot S_A)\Gamma(a)$ in this proof. In the case of a comodule {\it bi}-algebra, we find that this formula implies that $(1\ot S_A)\Gamma(a)=\Gamma(S_A(a))$ as expected. It is the dual version of the equality $S_B(b\tl a)=S_B(b)\tl a$ that we used in the proof of Proposition 1.3 for a module bi-algebra. As it can be seen above, in the general case, the formula is more complicated. Also this more general formula has a dual form that could be used to obtain a proof of the equality $\rho=\eta$ in Proposition 1.3, without the extra assumption that $B$ is an $A$-module bi-algebra. However, we will wait for later results on duality in order to give a completely {\it rigorous proof} of this result (see Corollary 3.8 in Section 3). 
\nl
We will now prove a few more properties of this {\it distinguished multiplier} $y$ satisfying and characterized by
the formulas 
$(\iota\ot\varphi)\Gamma(a)=\varphi(a)y$ and $(\iota\ot\psi)\Gamma(a)=\psi(a)y$ 
for all $a\in A$.
\snl
The first one is easy.

\inspr{1.6} Proposition \rm
The element $y$ is group-like. So we have $\Delta_B(y)=y\ot y$, as well as $\varepsilon_B(y)=1$ and $S_B(y)=y^{-1}$.
\snl \bf Proof\rm:
Take any $a\in A$ and apply $\varphi$ to the third leg of the equality
$$(\Delta_B\ot\iota_A)\Gamma(a)=(\iota_B\ot \Gamma)\Gamma(a).$$
Remember that this formula is true because we have a left coaction (see Definition 1.6 in [De-VD-W]). 
This gives
$$\align \varphi(a)\Delta_B(y) &= (\iota_B \ot \iota_B \ot \varphi)((\iota_B\ot \Gamma)\Gamma(a)) \\
   &= ((\iota_B\ot\varphi) \Gamma(a)) \ot y \\
   &= \varphi(a) y\ot y.
\endalign$$
The other statements follow easily because $y$ is group-like.
\hfill $\square$
\einspr

The next result is less obvious. We will show that $y$ commutes with the modular element $\delta_A$ of $A$ in the multiplier algebra $M(AB)$. 

\inspr{1.7} Proposition \rm
For all $b\in B$ we have
$$y^{-1}b y = b\tl \delta_A = \delta_A^{-1}b\delta_A.$$

\snl \bf Proof\rm:
Start with the equality
$$\sum_{(a)} (b\tl a_{(1)}\ot 1)
\Gamma (a_{(2)})= \sum_{(a)} 
\Gamma(a_{(1)})(b\tl a_{(2)}\ot 1),$$
true for all $a\in A$ and $b\in B$. It is one of the basic assumptions for a matched pair (see Section 2, in particular Assumption 2.12 and Proposition 2.3 in [De-VD-W]).
\snl
If we apply $\iota\ot \varphi$ we find
$$\varphi(a)(b\tl 1)y=\varphi(a) y (b\tl \delta_A).$$
So $y^{-1}b y = b\tl \delta_A$. Because $\delta_A$ is group-like, we have also $b\delta_A=\delta_A(b\tl \delta_A)$ and so $b\tl \delta_A = \delta_A^{-1}b\delta_A$ in $AB$.
\hfill $\square$\einspr

If we replace $b$ by $y$, we get that $\delta_A$ and $y$ commute with each other. 
\snl
Because $\delta_A$ is group-like in $A$, we should have $\Delta_\#(\delta_A)=(\delta_A\ot 1)
\Gamma(\delta_A)$. This last formula can be used to define $\Gamma(\delta_A)$. Observe that we can not extend $\Gamma$ to $M(A)$ in an obvious way.
\snl
Remember that $\delta_A$ is group-like for the coproduct $\Delta_A$ on $A$ but that the coproduct $\Delta_\#$ does not coincide with $\Delta_A$ on $A$. Still we have the following.

\inspr{1.8} Proposition \rm
We have $\Delta_\#(\delta_A)=\delta_A\ot \delta_A$.
\snl \bf Proof\rm:
Take $a\in A$. First, because $\varphi$ is a left integral on $A$, we get using the formula in Proposition 1.4 that
$$(\iota\ot \varphi)\Delta_\#(a\delta_A)=\sum_{(a)}a_{(1)}\delta_A\varphi(a_{(2)}\delta_A)y=\varphi(a\delta_A) y.$$
Similarly, now using that $\varphi(\,\cdot\,\delta_A)$ is a right integral and the formula in Proposition 1.5 we obtain
$$(\iota\ot \varphi)(\Delta_\#(a)(1\ot\delta_A))=\sum_{(a)}a_{(1)}\varphi(a_{(2)}\delta_A)y=\varphi(a\delta_A) \delta_A^{-1}y.$$
This implies that 
$$(\iota\ot \varphi)(\Delta_\#(a)(\delta_A\ot\delta_A))=\varphi(a\delta_A) y$$
and we get by comparing the two results that
$$(\iota\ot \varphi)\Delta_\#(a\delta_A)=(\iota\ot \varphi)(\Delta_\#(a)(\delta_A\ot\delta_A))$$
for all $a\in A$.
Now, it is not hard to see that the map $x\ot a\mapsto (x\ot 1)\Delta_\#(a)$ is bijective from $AB\ot A$ to itself. Therefore, multiplying the above equation with elements in $AB$, we get
$$(\iota\ot\varphi)((1\ot a)\Delta_\#(\delta_A))=\varphi(a\delta_A)\delta_A$$
for all $a$. Finally from the faithfulness of $\varphi$, we arrive at
$$\Delta_\#(\delta_A)=\delta_A \ot \delta_A.$$
In the proof, we have used that $y$ and $\delta_A$ commute.
\hfill $\square$
\einspr

So we see that $\delta_A$ is group-like, not only in $(A,\Delta_A)$, but also in $(AB,\Delta_\#)$. This is essentially the same as saying that we can write $
\Gamma(\delta_A)=1 \ot \delta_A$.
\nl
At the {\it end of this section}, we will give an example to illustrate already some of the results obtained here. Later, in Section 4, we will consider {\it more examples}.
\nl
{\it The $^*$-algebra case}
\nl
Now we assume that $A$ and $B$ are multiplier Hopf $^*$-algebras. For the right action $\tl$ of $A$ on $B$ we have the {\it extra assumption}
$$(b\tl a)^*=b^* \tl S(a)^*$$
when $a\in A$ and $b\in B$. For the left coaction $
\Gamma$ of $B$ on $A$ we now {\it also require} 
$$
\Gamma(S(a)^*)=((\iota\ot S)
\Gamma(a))^*$$
when $a\in A$. We refer to Theorem 2.15 in [De-VD-W] for these formulas.
\snl
Look first at {\it the case where cointegrals exist}. We know that $h^*$ will be a right cointegral if $h$ is a left cointegral. Therefore, if we apply the involution to the formulas in Proposition 1.1, we will find
$$
\rho(S(a)^*)^-=\eta(a)$$
for all $a\in A$. This means that $
\rho^*=\eta$ when we define the involution on the dual of $A$ in the standard way. Combined with 
$\rho=\eta$, we would find that $
\rho$ is self-adjoint.
\snl
Next, take {\it the general case} again. We then get the following.

\inspr{1.9} Proposition \rm The element $y$, as obtained in Proposition 1.4, is self-adjoint.
\snl \bf Proof\rm:
Consider the left integral $\varphi$ and define $\overline\varphi$ by
 $\overline\varphi(a)=\varphi(a^*)^-$ for all $a\in A$. Apply the first formula of Proposition 1.4 to $S(a)^*$ and use that $\overline\varphi\circ S$ is a right integral. We get
$$\align \varphi(S(a)^*)y
       &=(\iota\ot\varphi) \Gamma(S(a)^*) \\
       &=(\iota\ot\varphi)(((\iota\ot S) \Gamma(a))^*)\\
       &=((\iota\ot\overline\varphi)((\iota\ot S)\Gamma(a)))^*\\
       &=\overline\varphi(S(a))^-y^*=(\varphi(S(a)^*)y^*.
\endalign$$
proving the result.
\hfill $\square$\einspr

One can verify that the self-adjointness of $y$ is in agreement with the equation obtained in Proposition 1.7.
\nl
We finally look at some {\it examples}.

\iinspr{1.10} Example \rm 
i) Take a regular multiplier Hopf algebra $(A,\Delta_A)$ and let $(B,\Delta_B)=(A,\Delta_A^{\text{cop}})$. Define the right action $\tl$ of $A$ on $B$ by 
$$b\tl a=\sum_{(a)} S_A(a_{(1)})ba_{(2)}.$$
For this example, we will use $A\# B$ to denote the smash product (and {\it not} $AB$ for obvious reasons).
The left coaction $\Gamma$ of $B$ on $A$ is defined by
$$\Gamma(a)=\sum_{(a)} S_A(a_{(1)})a_{(3)} \ot a_{(2)}.$$
{\it We use Sweedler notation (only) for} $\Delta_A$. 
\snl
One can verify that, just as in the case of Hopf algebras (see e.g.\ Example 6.2.8 in [M]), this gives a matched pair. It is not hard to prove that the action is unital and that $B$ is a right $A$-module algebra. Similarly, it is relative easy to show that $\Gamma$ makes $A$ into a left $B$-comodule coalgebra. It is more complicated in this case to show that the action and the coaction are compatible so that the coproduct $\Delta_\#$ on $A\# B$ is an algebra map. We will consider this case again in Section 4 where we study more examples and see that there is a better way to show this last property (see a remark following the proof of Proposition 4.2).
\snl
If we started above with a multiplier Hopf $^*$-algebra, one easily verifies that $(b\tl a)^*=b^*\tl S_A(a)^*$ for all $a\in A$ and $b\in B$. Also $\Gamma(S_A(a)^*)=((\iota\ot S_A)\Gamma(a))^*$ for all $a$. These are the compatibility conditions with the involutive structure.
\snl
ii) Now suppose that we have integrals on $A$. Let $\varphi_A$ and $\psi_A$ be a left and a right integral. We find
$$\align (\iota\ot\varphi_A)\Gamma(a) 
      &= \sum_{(a)}\varphi_A(a_{(2)})S_A(a_{(1)})a_{(3)}\\
      &= \sum_{(a)}\varphi_A(a_{(1)})a_{(2)}\\
		  &=\varphi_A(a)\delta\\
  (\iota \ot\psi_A)\Gamma(a)
      &= \sum_{(a)}\psi_A(a_{(2)})S_A(a_{(1)})a_{(3)}\\
      &= \sum_{(a)}\psi_A(a_{(2)})S_A(a_{(1)})\\
      &=\psi_A(a)S_A(\delta^{-1})=\psi_A(a)\delta.
\endalign$$
In this calculation, $\delta=\delta_A$. 
\snl
This nicely illustrates the property $y=z$. Here we find that this distinguished multiplier in $M(A)$ coincides with the modular element $\delta$ as expected. Observe that indeed $b\tl\delta=S(\delta)b\delta=\delta^{-1}b\delta$. The formula can be interpreted, not only in $A$, but also in $A\# B$ because
$$(\delta^{-1}\# 1)(1\# b)(\delta\# 1)=(\delta^{-1}\# b)(\delta\# 1)=\delta^{-1}\delta \#(b\tl \delta)=1\#\delta^{-1}b\delta.$$
Remark finally that (at least formally)
$$\Gamma(\delta)=S(\delta)\delta\ot\delta=\delta^{-1}\delta\ot \delta=1\ot \delta.$$
\vskip -0.7 cm\hfill$\square$\einspr

One might also consider the {\it dual example}. Then we start with a regular multiplier Hopf algebra $(C,\Delta_C)$ and we let $(D,\Delta_D)=(C^{\text{op}},\Delta_C)$. The left action of $D$ on $C$ is defined as before by the formula
$d\tr c=\sum_{(d)} S_C(d_{(1)})cd_{(2)}$
and the right coaction of $C$ on $D$ by 
$$\Gamma(d)=\sum_{(d)}d_{(2)}\ot S_C(d_{(1)})d_{(3)}.$$
Now suppose that we have integrals on $C$. Let $\varphi$ and $\psi$ be a left and a right integral on $C$. They are also left and right integrals on $D$ because $\Delta_D=\Delta_C$ is assumed. We easily calculate that $(\varphi\ot\iota)\Gamma(d)=\varphi(d)\delta$ as well as $(\psi\ot\iota)\Gamma(d)=\psi(d)\delta$ where $\delta=\delta_C=\delta_D$.
\snl
This example can either be obtained from the other one by duality or by the conversion procedure as explained in [De-VD-W] (see a remark in Section 1 of that paper).
We will also consider this example again in Section 4.
\nl\nl

\bf 2. Integrals on the bicrossproduct (and related objects) \rm
\nl
As in the previous section, $A$ and $B$ are regular multiplier Hopf algebras with a right action of $A$ on $B$ and a left coaction of $B$ on $A$ making $(A,B)$ into a matched pair. In this section, we will moreover {\it assume that both $A$ and $B$ have integrals}, i.e.\ that they are algebraic quantum groups. We have shown already in [De-VD-W] that then also the bicrossproduct has integrals. The formula for the right integral was given in Theorem 3.3 of [De-VD-W]. In this section, we will first recall this result. We will however also give the formula for the left integral. And we will find formulas for the data associated with these integrals on the bicrossproduct.
\snl
We denote a left integral by $\varphi$ and a right integral by $\psi$. So, we use $\varphi_A, \varphi_B$ for left integrals on $A$ and $B$ respectively and similarly  $\psi_A, \psi_B$ for right integrals. We will make {\it the convention} that $\varphi=\psi\circ S$, both for $A$ and for $B$. Observe that this is different from the rather common convention where it is assumed that $\psi=\varphi\circ S$. We also use $\delta$ for the modular element and $\sigma$ and $\sigma'$ for the modular automorphisms associated with $\varphi$ and $\psi$ respectively, again with the appropriate indices if necessary. For the corresponding objects on the smash product $AB$, we write $\varphi_\#, \psi_\#, \dots$.
\nl
The formula for {\it the right integral} $\psi_\#$ on the bicrossproduct was already obtained in Theorem 3.3 of [De-VD-W]. We recall it here for completeness.

\inspr{2.1} Proposition \rm
If we define $\psi_\#$ on $AB$ by $\psi_\#(ab)=\psi_A(a)\psi_B(b)$, then we have a right integral on the bicrossproduct $(AB,\Delta_\#)$. 
\hfill $\square$\einspr

Formally, the argument goes as follows. Take $a\in A$ and $b\in B$. Then
$$\align (\psi_\#\ot\iota)\Delta_\#(ab)
   &=\sum_{(a)} (\psi_\#\ot \iota)((a_{(1)} \ot 1)\Gamma(a_{(2)})\Delta_B(b)) \\
   &=\sum_{(a)} \psi_A(a_{(1)})(\psi_B \ot \iota)(\Gamma(a_{(2)})\Delta_B(b)) \\
   &=\psi_A(a)(\psi_B\ot\iota)(\Gamma(1)\Delta_B(b)) \\
   &=\psi_A(a)(\psi_B\ot\iota)\Delta_B(b) \\
   &=\psi_A(a)\psi_B(b)1.
\endalign$$
Of course, the problem is that $\Gamma(1)$ is not really defined and moreover it is unclear if we can use that 
$\sum_{(a)} \psi_A(a_{(1)})\Gamma(a_{(2)}) =\psi_A(a)\Gamma(1)$. In the proof of Theorem 3.3 in [De-VD-W] it is shown how this problem can be overcome. The basic idea is to use as a starting point not the previous formula but $\sum_{(a)} \psi_A(a_{(1)})\Delta_\#(a'a_{(2)}) =\psi_A(a)\Delta_\#(a')$ as this makes perfectly good sense for all $a,a'\in A$. Moreover, although we can not talk about $\Gamma(1)$, we can say that $\Delta_\#(1)=1$ and this should mean precisely that $\Gamma(1)=1$.
\nl
There are a few easy properties of the right integral $\psi_\#$ obtained above.

\inspr{2.2} Proposition \rm
For $a\in A$ and $b\in B$ we also have $\psi_\#(ba)=\psi_A(a)\psi_B(b)$. In fact, we find $\psi_\#(bab')=\psi_A(a)\psi_B(bb')$ for all $a\in A$ and $b,b'\in B$.
\hfill $\square$\einspr

Indeed, as $ba=\sum_{(a)}a_{(1)}(b\tl a_{(2)})$ we get
$$\align \psi_\#(bab')
   &=\sum_{(a)}\psi_A(a_{(1)})\psi_B((b\tl a_{(2)})b') \\
   &=\psi_A(a)\psi_B((b\tl 1)b')=\psi_A(a)\psi_B(bb').
\endalign$$
\nl
We now look for {\it the left integral} $\varphi_\#$ on the bicrossproduct $(AB,\Delta_\#)$ and {\it the modular element} $\delta_\#$. 

\inspr{2.3} Proposition \rm
If we define the left integral $\varphi_\#$ as $\psi_\#\circ S_\#$, it is given by the formula
$$\varphi_\#(ab)=\varphi_A(a)\varphi_B(yb)$$
where $y$ is the element in $M(B)$, defined in Proposition 1.4 by the equation \newline 
$(\iota\ot\varphi_A)\Gamma(a)=\varphi_A(a)y$ for all $a\in A$. The {\it modular element} is given by $\delta_\#=\delta_B\delta_A y^{-1}$.
\snl \bf Proof\rm:
As already mentioned before, we let $\varphi_\#=\psi_\#\circ S_\#$. We now use the formula for $S_\#$ as obtained in Theorem 3.1 of [De-VD-W], recalling also the convention for the Sweedler type notation $\Gamma(a)=\sum_{(a)}a_{(-1)}\ot a_{(0)}$ (see also Section 1 of [De-VD-W]). We find
$$\align \varphi_\#(ab)
     &=\sum_{(a)} \psi_\#(S_B(a_{(-1)}b)S_A(a_{(0)})) \\
     &=\sum_{(a)} \psi_B(S_B(a_{(-1)}b))\psi_A(S_A(a_{(0)})) \\
     &=\psi_B(S_B(yb))\psi_A(S_A(a))\\
     &=\varphi_B(yb)\varphi_A(a)
\endalign$$
for all $a\in A$ and $b\in B$. We have used that $\varphi_A=\psi_A\circ S_A$ and similarly for $B$.
\snl
In order to find the modular element, we use that $(\iota\ot \psi_\#)\Delta_\#(ab)=\psi_\#(ab)\delta_\#^{-1}$ whenever $a\in A$ and $b\in B$. We get
$$\align (\iota\ot\psi_\#)\Delta_\#(ab)
     &=\sum_{(a)} (\iota\ot\psi_\#)((a_{(1)}\ot 1)\Gamma(a_{(2)})\Delta_B(b)) \\
     &=\sum_{(a)} a_{(1)}((\iota\ot\psi_A)\Gamma(a_{(2)}))((\iota\ot\psi_B)\Delta_B(b)) \\
     &=\sum_{(a)} a_{(1)}\psi_A(a_{(2)})y((\iota\ot\psi_B)\Delta_B(b)) \\
     &=\psi_A(a)\psi_B(b)\delta_A^{-1}y\delta_B^{-1}
\endalign$$
so that $\delta_\#=\delta_B y^{-1} \delta_A$. As we know that $\delta_A$ and $y$ commute (see a remark following Proposition 1.7), we get the desired result.
\hfill $\square$\einspr

It is instructive to {\it verify some formulas}. The modular element $\delta_\#$  should also verify e.g.\ the equation $(\varphi_\#\ot \iota)\Delta_\#(ab)=\varphi_\#(ab)\delta_\#$ for all $a\in A$ and $b\in B$. To see this, take $a\in A$ and $b\in B$. We find (at least formally),
$$\align (\varphi_\#\ot \iota)\Delta_\#(ab)
     &=\sum_{(a)} (\varphi_\#\ot \iota)((a_{(1)}\ot 1)\Gamma(a_{(2)})\Delta_B(b)) \\
     &=\sum_{(a)} \varphi_A(a_{(1)})(\varphi_B\ot \iota)((y\ot 1)\Gamma(a_{(2)})\Delta_B(b)) \\
     &=\varphi_A(a)(\varphi_B\ot \iota)((y\ot 1)\Gamma(\delta_A)\Delta(b)) \\
     &=\varphi_A(a)\delta_A (\varphi_B\ot \iota)((y\ot 1)\Delta(b)) \\
     &=\varphi_A(a)\delta_A \varphi_B(yb)y^{-1} \delta_B.
\endalign$$
Finally, we know  that $\delta_A y^{-1}$ commutes with all elements of $B$ (again from Proposition 1.7) and that also $y$ commutes with $\delta_A$. So 
$$\delta_A y^{-1} \delta_B=\delta_B y^{-1}\delta_A.$$
We see again that $\delta_\#=\delta_B y^{-1}\delta_A$ as asserted.
\snl
In the above formal argument we have used $\Gamma(\delta_A)=1\ot \delta_A$, together with the fact that $\sum_{(a)}\varphi_A(a_{(1)})\Gamma(a_{(2)})=\varphi_A(a)\Gamma(\delta_A)$. The problem here is the same as with the formal proof of the right invariance of $\psi_\#$ as explained in the beginning of this section. Also here the argument can be made precise if in stead we use $\sum_{(a)}\varphi_A(a_{(1)})\Delta_\#(a'a_{(2)})=\varphi_A(a)\Delta_\#(a'\delta_A)$, but there is no need to explain this argument in more detail. We also have used that $y$ is group-like in $M(B)$ (see Proposition 1.6).
\snl

We  want to include the following remark.

\inspr{2.4} Remark \rm
i) In an earlier version of the paper on bicrossproducts for multiplier Hopf algebras (unpublished), a formula for $\delta_\#$ had been obtained under an extra assumption. It was shown that $\delta_\#=\delta_A(\delta_B\tl \delta_A)$ provided
$$(\iota\ot\psi_A)(\Gamma(a)(b\ot 1))=\psi_A(a)(b\tl \delta_A)$$
is assumed for all $a\in A$ and $b\in B$. Using our result, this condition simply would mean that $yb=b\tl \delta_A$ and because we have seen in our Proposition 1.7 that $b\tl \delta_A=y^{-1}by$, the condition would imply that $y=1$. Then the formula $\delta_\#=\delta_A(\delta_B\tl \delta_A)$ gives $\delta_\#=\delta_A\delta_B$. We see from this that the condition is rather artificial (and not needed) for this result.
\snl
ii) Also in Proposition 5.4 of [B-*], the modular element $\delta_\#$ is found in the case where $A$ and $B$ are Hopf algebras. However, the formula in that paper is also not natural because the multiplier $y$ we obtained here in Proposition 1.4 is missing.
\hfill$\square$\einspr

We now get the {\it scaling constant} $\tau_\#$.

\inspr{2.5} Proposition \rm
We have $\tau_\#=\tau_A\tau_B$ where $\tau_A$ and $\tau_B$ are the scaling constants for $A$ and $B$ respectively.
\snl\bf Proof\rm:
Although the result looks quite obvious and simple, this is not true for the proof of it. We use the formula that defines the scaling constant, namely
$\psi_\#(S_\#^2(x))=\tau_\#\psi_\#(x)$ when $x\in AB$. The difficulty with this is the complexity of the antipode $S_\#$. 
\snl
So let $x=ab$ with $a\in A$ and $b\in B$. Then we get, using the formula for $S_\#$, 
$$\align \psi_\#(S_\#^2(ab))
   &=\sum_{(a)}\psi_\#(S_\#(S_B(a_{(-1)}b)S_A(a_{(0)})))\\
   &=\sum_{(a)}\psi_\#(S_\#(S_A(a_{(0)}))S_B^2(a_{(-1)}b))\\
   &=\sum_{(a)}\psi_\#(S_B(S_A(a_{(0)})_{(-1)})S_A(S_A(a_{(0)})_{(0)})S^2_B(a_{(-1)}b)).
\endalign$$
Using  that $\psi_\#(b'ab)=\psi_A(a)\psi_B(b'b)$ for all $a\in A$ and $b,b'\in B$ and also that 
$$(\iota \ot (\psi_A\circ S_A))\Gamma(a')=(\psi_A\circ S_A)(a')y$$
for all $a'\in A$, we find
$$\align \psi_\#(S_\#^2(ab))
    &=\sum_{(a)}\psi_B(S_B(S_A(a_{(0)})_{(-1)})S^2_B(a_{(-1)}b))\psi_A(S_A(S_A(a_{(0)})_{(0)})) \\
    &=\sum_{(a)}\psi_B(S_B(y)S^2_B(a_{(-1)}b))\psi_A(S_A^2(a_{(0)}))\\
    &=\psi_B(S_B(y)S^2_B(yb))\psi_A(S_A^2(a))\\
    &=\psi_B(S^2_B(b))\psi_A(S_A^2(a))\\
    &=\tau_B\tau_A\psi_B(b)\psi_A(a)=\tau_B\tau_A\psi_\#(ab).
\endalign$$ 
This proves the result.
\hfill $\square$
\einspr

The reader can verify the necessary coverings. 
\snl
There are other possibilities to obtain this result. One such possibility would be to use that $\sigma'_\#(\delta_\#)=\tau_\#^{-1} \delta_\#$. Unfortunately, for this we need the formula for $\sigma'_\#(\delta_A)$ and that will only be obtained later (see a remark following Proposition 2.10). 
\nl
Finally, we try to find {\it formulas for the modular automorphisms} $\sigma_\#$ and $\sigma'_\#$. Because of the formulas in the Propositions 2.1, 2.2 and 2.3, it is relatively easy to see what the modular automorphisms do on elements of $B$. This is the content of the next proposition.

\inspr{2.6} Proposition \rm For $b\in B$  we have
$$\sigma_\#(b)=\sigma_B(b) \qquad\qquad \text{and} \qquad\qquad \sigma'_\#(b)=\sigma'_B(b).$$
\bf\snl Proof\rm:
We first prove the last formula. We have for all $a'\in A$  and $b,b'\in B$ that
$$\align \psi_\#(a'b'\sigma'_\#(b))&=\psi_\#(ba'b')=\psi_A(a')\psi_B(bb') \\
             &=\psi_A(a')\psi_B(b'\sigma'_B(b))=\psi_\#(a'b'\sigma'_B(b)).
\endalign$$
To obtain the other formula, we use that $\sigma_\#(b)=\delta_\#^{-1}\sigma'_\#(b)\delta_\#$ for all $b$. We get using the formula for $\delta_\#$, given in Proposition 2.3 that
$$\align \sigma_\#(b)&=\delta_\#^{-1}\sigma'_\#(b)\delta_\# \\
    &=\delta_A^{-1}y\delta_B^{-1}\sigma'_B(b)\delta_By^{-1}\delta_A\\
    &=\delta_A^{-1}y\sigma_B(b)y^{-1}\delta_A
\endalign$$
and the result follows as we have seen that $y^{-1}\delta_A$ commutes with $B$.
\hfill$\square$\einspr

One can also obtain the result directly, using the formula for $\varphi_\#$, obtained in Proposition 2.3, in combination with the result of Proposition 1.7.
\snl
The formulas for $\sigma_\#(a)$ and $\sigma'_\#(a)$ with $a\in A$ are more complicated. 
\snl
We first prove two lemmas. The results will be used also later in Section 3.

\inspr{2.7} Lemma \rm Given $a'\in A$ and $b'\in B$, there exists $b''\in B$ so that 
$$\psi_B((b\tl a')b')=\psi_B(bb'')$$
for all $b\in B$.
\snl\bf Proof\rm:
To prove this result, take another element $p\in A$ and assume that $\psi_A(p)=1$. Then, for all $b\in B$, we have
$$\align \psi_B((b\tl a')b')
	&=\psi_A(p)\psi_B((b\tl a')b')\\
	&=\psi_\#(p(b\tl a')b')\\
	&=\sum_{(a')}\psi_\#(pS_A(a'_{(1)})ba'_{(2)}b')\\
	&=\sum_{(a')}\psi_\#(ba'_{(2)}b'\sigma'_\#(pS_A(a'_{(1)}))).\\
\endalign$$
We  see that there exists an element $z\in AB$  so that $\psi_B((b\tl a')b')=\psi_\#(bz)$ for all $b\in B$. If we write $z=\sum_i b_ia_i$ with $a_i\in A$ and $b_i\in B$, we find that
$$\psi_B((b\tl a')b')=\sum_i \psi_A(a_i)\psi_B(b b_i)$$
for all $b\in B$. If we let $b''=\sum_i\psi_A(a_i)b_i$, we have proven the formula in the lemma.

\hfill$\square$\einspr

The above result is the same as saying that $b\mapsto \omega(b\tl a)$ belongs to the dual $\widehat B$ for all $\omega\in\widehat B$ and $a\in A$. In fact, from the proof we can already see that all elements in $\widehat B$ can be obtained as linear combinations of such linear functionals $\omega((\,\cdot\,)\tl a)$ with $\omega\in\widehat B$ and $a\in A$. We will give a rigorous proof of this statement in the  next section (see Theorem 3.3).
\snl
As a consequence, we also get the following, similar properties.

\inspr{2.8} Lemma \rm
Given $a\in A$ and $b'_1,b'_2\in B$, there exist elements $b_1'',b_2''\in B$ so that 
$$\align \psi_B((bb'_1)\tl a)&=\psi_B(bb_1'')\\
         \psi_B((b'_2b)\tl a)&=\psi_B(bb_2'')
\endalign$$
for all $b\in B$.
\snl\bf Proof\rm:
Indeed, in the first case we use
$$\psi_B((bb'_1)\tl a)=\sum_{(a)}\psi_B((b\tl a_{(1)})(b'_1\tl a_{(2)}))$$ 
and apply the previous lemma with $a_{(1)}$ and $b'_1\tl a_{(2)}$ in the place of $a'$ and $b'$. For the second case, we use
$$\align \psi_B((b'_2b)\tl a)
            &=\sum_{(a)}\psi_B((b'_2\tl a_{(1)})(b\tl a_{(2)}))\\
            &=\sum_{(a)}\psi_B((b\tl a_{(2)})\sigma'_B(b'_2\tl a_{(1)}))
\endalign$$
and again apply the previous lemma.
\hfill$\square$\einspr

These results say that $\psi_B((\,\cdot\,b)\tl a)$ and $\psi_B((b\,\cdot\,)\tl a)$ are elements of $\widehat B$ for all $a\in A$ and $b\in B$.
\snl
Now, from the very existence of $\sigma'_\#(a)$, we  obtain the following result. It will be used to give the  formula for  $\sigma'_\#(a)$ in Proposition 2.10 below.

\inspr{2.9} Proposition \rm There exists a linear map $\gamma$ from $A$ to $M(B)$ such that $\psi_B(b\tl a)=\psi_B(b\gamma(a))$ for all $a\in A$ and $b\in B$.

\bf\snl Proof\rm: 
Given $a\in A$, we define the multiplier $\gamma(a)\in M(B)$, using the notations as above, in Lemma 2.8, by
$$b'_1\gamma(a)=b_1''\qquad\qquad \text{and}\qquad\qquad \gamma(a)\sigma'_B(b'_2)=b_2''.$$
A straightforward calculation gives that 
$$ (b'_1\gamma(a))\sigma'_B(b'_2)=b'_1(\gamma(a)\sigma'_B(b'_2))$$
and this proves that $\gamma(a)$ is indeed well-defined as an element in $M(B)$. 
\hfill$\square$\einspr

Now, we can give the formula for $\sigma'_\#(a)$ when $a\in A$.

\iinspr{2.10} Proposition \rm
For all $a\in A$  we have
$$\sigma'_\#(a)=\sum_{(a)} \sigma'_A(a_{(1)})\gamma(S^{-1}(a_{(2)})).$$
\bf\snl Proof\rm:
Take $a,a'\in A$ and $b\in B$. Then we have
$$\align \psi_\#(aba')
      &=\sum_{(a')}\psi_\#(aa'_{(1)}(b\tl a'_{(2)}))\\
      &=\sum_{(a')}\psi_A(aa'_{(1)})\psi_B(b\tl a'_{(2)})\\
      &=\sum_{(a')}\psi_A(aa'_{(1)})\psi_B(b\gamma(a'_{(2)})).\endalign$$
Now we use the formula
$$\sum_{(a')}\psi_A(aa'_{(1)})a'_{(2)}=\sum_{(a)}\psi_A(a_{(1)}a')S^{-1}(a_{(2)}).$$
Then we get
$$\align \psi_\#(aba')
      &=\sum_{(a)}\psi_A(a_{(1)}a')\psi_B(b\gamma(S^{-1}(a_{(2)})))\\
      &=\sum_{(a)}\psi_A(a'\sigma'_A(a_{(1)}))\psi_B(b\gamma(S^{-1}(a_{(2)})))\\
      &=\sum_{(a)}\psi_\#(ba'\sigma'_A(a_{(1)})\gamma(S^{-1}(a_{(2)})))
\endalign$$
and this proves the result.
\hfill$\square$
\einspr

For $\sigma_\#(a)$ we get $\delta_\#^{-1}\sigma'_\#(a)\delta_\#$ and if we insert the given formulas for $\delta_\#$ and $\sigma'_\#(a)$, we get the one for $\sigma_\#(a)$. However, we do not seem to get any simpler expressions as we do in the case of $\sigma_\#$ on $B$.
\snl
As we promised before, let us consider (and verify formally) the equation $\sigma'_\#(\delta_\#)=\tau_\#^{-1} \delta_\#$. We have
$$\sigma'_\#(\delta_\#)=\sigma'_\#(\delta_Ay^{-1}\delta_B)=\sigma'_\#(\delta_A)\sigma'_B(y^{-1}\delta_B).$$
From the previous proposition, we know that $\sigma'_\#(\delta_A)=\sigma'_A(\delta_A)\gamma(\delta^{-1}_A)$. From the definition of $\gamma$, we get
$$\psi_B(b\gamma(\delta_A^{-1}))=\psi_B(b\tl\delta_A^{-1})=\psi_B(yby^{-1})=\psi_B(by^{-1}\sigma'_B(y))$$
for all $b\in B$ so that $\gamma(\delta_A^{-1})=y^{-1}\sigma'_B(y)$. Combining all these results we find
$$\align \sigma'_\#(\delta_\#)
                   &=\sigma'_\#(\delta_A)\sigma'_B(y^{-1})\sigma'_B(\delta_B)\\
                   &=\sigma'_A(\delta_A)y^{-1}\sigma'_B(y)\sigma'_B(y^{-1})\sigma'_B(\delta_B)\\
                   &=\tau_A^{-1}\delta_A y^{-1}\tau_B^{-1}(\delta_B)
\endalign$$
and this is precisely what we want.
\nl
Before we continue, we need to give a couple of remarks on the results obtained in the Propositions 2.9 and 2.10.

\iinspr{2.11} Remark \rm
i) If the right integral $\psi_B$ is {\it invariant under the action} of $A$, that is if 
$\psi_B(b\tl a)=\varepsilon_A(a)\psi_B(b)$
for all $a\in A$ and $b\in B$, we see that $\gamma(a)=\varepsilon_A(a)1$ and we find that also
$\sigma'_\#(a)=\sigma'_A(a)$ for all $a\in A$. A simple calculation gives that $\sigma_\#(a)=\delta_B^{-1}y\sigma_A(a)y^{-1}\delta_B$ and so in general we need not have that also $\sigma_\#$ concides with $\sigma_A$ on $A$ or that it even maps $A$ into itself. 
\snl
ii) If $B$ is a right $A$-module {\it bi}-algebra, that is if 
$$\Delta_B(b\tl a)= \sum_{(a)(b)}(b_{(1)}\tl a_{(1)})\ot (b_{(2)}\tl a_{(2)}),$$
one can verify that $b\mapsto \psi_B(b\tl a)$ is also right invariant, and so of the form $\gamma_0(a)\psi_B$ where now $\gamma_0$ is a homomorphism of $A$ to $\Bbb C$. Then $\gamma(a)=\gamma_0(a)1$ and we see that
$$\sigma'_\#(a)=\sum_{(a)}\gamma_0(S^{-1}(a_{(2)}))\sigma'_A(a_{(1)}).$$
In this case, still $\sigma'_\#$ will leave $A$ invariant but it no longer coincides with $\sigma'_A$.
\hfill$\square$
\einspr

At the end of this section, we will also look at the map $\gamma$ for the example with the adjoint action and the adjoint coaction (Example 1.10 in the previous section). 
\snl
Before, let us look at {\it the $^*$-algebra case}. Not much can be said except for the following important result.

\iinspr{2.12} Proposition \rm If $A$ and $B$ are multiplier Hopf $^*$-algebras with positive integrals, then also the bicrossproduct has positive integrals.

\snl\bf Proof: \rm
For all $a,a'\in A$ and $b,b'\in B$ we have, using an earlier result,
$$\psi_\#((ab)^*(a'b'))=\psi_\#(b^*a^*a'b')=\psi_A(a^*a')\psi_B(b^*b')$$
and then the result follows.
\hfill $\square$
\einspr

We know that in this case, the scaling constants are $1$ (see e.g.\ Theorem 3.4 in [DC-VD]). This is in accordance with the formula of Proposition 2.5. 
\snl
Also, we know that the modular elements are self-adjoint. We have seen that $\delta_\#=\delta_By^{-1}\delta_A$ and so self-ajointness of $\delta_\#$ is equivalent with the equality 
$$\delta_By^{-1}\delta_A=\delta_Ay^{-1}\delta_B$$
because we know that also $y$ is self-adjoint (cf.\ Proposition 1.9). We have seen already that the above equality is satisfied. See some remarks following Proposition 2.3.
\nl
Now we finish with {\it the examples}. We use the notations and conventions of Example 1.10. In particular, we will not use $AB$ but $A\# B$ for the smash product. And we use the Sweedler notation for $\Delta_A$.

\iinspr{2.13} Example \rm 
i) First observe that $\psi_B=\varphi_A$ and $\varphi_B=\psi_A$ because $(B,\Delta_B)=(A,\Delta_A^{\text{cop}})$. Then we find that $\psi_\#(a\# b)=\psi_A(a)\varphi_A(b)$ for $a,b\in A$. As we know that the distinguished multiplier $y\in M(B)$ in this case is $\delta_A$, we find that $\varphi_\#(a\# b)=\varphi_A(a)\psi_A(\delta_A b)$.
\snl
ii) For the modular element $\delta_\#$, we use the formula $\delta_\#=\delta_A \# (y^{-1}\delta_B)$. And because $\delta_B=\delta^{-1}_A$, we get $\delta_\#=\delta_A \# \delta^{-2}_A$.
\snl
iii) We now look for the modular automorphisms. First, we calculate and discuss the map $\gamma: A\to M(B)$ as obtained in Proposition 2.9. We will use $S$ for the antipode $S_A$ of $(A,\Delta_A)$. We find 
$$\psi_B(b\tl a)=\sum_{(a)}\varphi_A(S(a_{(1)})ba_{(2)})=\sum_{(a)}\varphi_A(b a_{(2)}\sigma_A(S(a_{(1)})))$$
 and we see that
$$\gamma(a)=\sum_{(a)}a_{(2)}\sigma_A(S(a_{(1)}))$$
for all $a\in A$. 
\snl
If we assume that the fourth power of $S$ equals $\iota$, we can rewrite this as
$$\align \gamma(a)
                    &=\sum_{(S(a))}(S^{-1}(S(a)_{(1)}))\sigma_A(S(a)_{(2)}) \\
                    &=\sum_{(S(a))}(S^{-3}S^2(S(a)_{(1)}))\sigma_A(S(a)_{(2)}) \\
                    &=\sum_{(\sigma_A(S(a)))}(S^{-3}(\sigma_A(S(a))_{(1)}))(\sigma_A(S(a))_{(2)}) \\
                    &=\sum_{(\sigma_A(S(a)))}(S(\sigma_A(S(a))_{(1)}))(\sigma_A(S(a))_{(2)}) \\
                    &=\varepsilon_A(\sigma_A(S(a)))1_A.
\endalign$$
We have used that in general $\Delta(\sigma(a))=(S^2\otimes \sigma)\Delta(a)$. Finally, because also in general $\sigma\circ S=S\circ {\sigma'}^{-1}$ we get
$$\gamma(a)=\varepsilon_A({\sigma_A'}^{-1}(a))1_A=\varepsilon_A(\sigma_A^{-1}(a))1_A.$$
\snl
One can see that also conversely, if $\gamma(a)=\varepsilon_A(\sigma_A^{-1}(a))1_A$ for all $a$, then we must have that $S^4=\iota_A$. Recall that $\varepsilon_A\circ\sigma_A^{-1}$ is the homomorphism on $A$ given by the pairing with the modular element $\delta_{\widehat A}$ of the dual $\widehat A$ of $A$. This gives an example of the case where $\psi_B$ will be {\it relatively invariant under the action} of $A$ and not necessarily invariant.
\snl
iv) Now, no longer assume this extra condition. Then, we get for the modular automorphism $\sigma'_\#$ the following:
$$\align \sigma'_\#(a\#b)
         &=\sum_{(a)} \sigma'_A(a_{(1)}) \# \gamma(S^{-1}(a_{(2)}))\sigma'_B(b) \\
         &=\sum_{(a)} \sigma'_A(a_{(1)}) \# S^{-1}(a_{(2)})\sigma_A(a_{(3)})\sigma'_B(b)
\endalign$$
for all $a,b\in A$. In the event that $S^4=\iota$, we find that $\sigma'_\#$ maps $A$ to $A$ and that
$\sigma'_\#=\sigma'_A\circ S^{-2}\circ\sigma_A$ on $A$. So, if moreover $\sigma_A=S^{-2}$ (as is the case for discrete quantum groups - see [VD3]), we will find that $\sigma'_\#(a)=\sigma_A(a)$ for all $a\in A$.
\hfill$\square$\einspr

Also here, we might consider the dual case where we start with a regular multiplier Hopf algebra 
$(C,\Delta_C)$ and where $(D,\Delta_D)=(C^{\text{op}},\Delta_C)$ as at the end of the previous section. 
\snl
We will again say {\it more about these two examples} in Section 4.
\nl\nl

\bf 3. The dual $(AB)\,{\widehat{}}$ of the bicrossproduct $AB$ \rm
\nl
In this section, we again take a matched pair $(A,B)$ of algebraic quantum groups as in the previous sections. The {\it main goal} is to show that the dual $(AB)\,{\widehat{}}\,$ of $AB$, in the sense of algebraic quantum groups, is again a bicrossproduct (of the second type) of the duals $\widehat A$ and $\widehat B$ (see Theorem 3.7 below). At the end of this section, we will also consider the $^*$-algebra case.
\snl
Conform with the notations used in [De-VD-W], we will denote $\widehat A$ also by $C$ and $\widehat B$ by $D$. For the various pairings we use $\langle\,\cdot\, , \,\cdot\,\rangle$. We will also use the conventions and notations about the integrals and the related objects for $A$, $B$ and $AB$ as explained in the beginning of the previous section.
\nl
{\it Formally}, the left coaction $\Gamma$ of $B$ on $A$ induces a left action of $D$ on $C$ by the formula
$$\langle a, d\tr c\rangle=\langle \Gamma(a),d\ot c\rangle$$
where $a\in A$, $c\in C$ and $d\in D$. Similarly, but again formally, the right action of $A$ on $B$ gives a right coaction $\Gamma$ of $C$ on $D$ by
$$\langle b\ot a,\Gamma(d)\rangle = \langle b\tl a,d \rangle$$
where $a\in A$, $b\in B$ and $d\in D$.
\snl
When $A$ and $B$ are finite-dimensional, these formulas will indeed make the pair $(C,D)$ into a matched pair (of the second type) and the bicrossproduct $CD$ of $C$ and $D$ is identified with the dual $(AB)\,{\widehat{}}\,$ by means of the natural pairing obtained from the pairing of $A$ with $\widehat A$ and $B$ with $\widehat B$. See e.g.\ Remark 1.7 in [De-VD-W].
\snl
However, in the more general (possibly non-finite-dimensional) case, even for ordinary Hopf algebras, the result is not so obvious. The problem is that the above formulas define $d\tr c\in A'$ and $\Gamma(d)\in (B\ot A)'$, the linear duals of $A$ and $B\ot A$ respectively, whereas we need these elements to belong to the (in general) strictly smaller spaces $\widehat A$ and $\widehat B\ot \widehat A$ respectively. 
\snl
We will focus in this section mostly on these problems and not so much on the (expected) formulas.
\nl
\it The identification of the dual $(AB)\,{\widehat{}}\,$ of $AB$ with $\widehat A\ot\widehat B$ as vector spaces\rm
\nl
First, we want to identify the dual $(AB)\,{\widehat{}}\,$ of $AB$ with the tensor product $\widehat A\ot\widehat B$ as linear spaces. This is the content of the following proposition.

\inspr{3.1} Proposition \rm Let $\omega'\in \widehat A$ and $\omega''\in \widehat B$ and define $\omega$ on $AB$ by $\omega(ab)=\omega'(a)\omega''(b)$. Then $\omega\in (AB)\,{\widehat{}}\,$ and the map $\omega'\ot\omega''\mapsto \omega$ gives a bijective linear map from $\widehat A\ot \widehat B$ to $(AB)\,{\widehat{}}\,$.
\snl\bf Proof\rm: 
Take elements $a,a'\in A$ and $b,b'\in B$ and let $\omega'=\psi_A(\,\cdot\,a')$ and $\omega''=\psi_B(\,\cdot\,b')$. Then
$$\align \omega'(a)\omega''(b)
			&=\psi_A(aa')\psi_B(bb') \\
      &=\psi_A({\sigma'_A}^{-1}(a')a)\psi_B(bb') \\
      &=\psi_\#({\sigma'_A}^{-1}(a')abb')\\
      &=\psi_\#(abb'\sigma'_\#({\sigma'_A}^{-1}(a'))) \\
      &=\psi_\#(ab\sigma'_\#({\sigma'_B}^{-1}(b'){\sigma'_A}^{-1}(a'))).
\endalign$$
This proves the first statement. Moreover, because the linear map from $A\ot B$ to $AB$, given by
$$a'\ot b'\mapsto \sigma'_\#({\sigma'_B}^{-1}(b'){\sigma'_A}^{-1}(a'))$$
is bijective, the result is proven.
\hfill $\square$  \einspr

There is also another possible identification of $(AB)\,{\widehat{}}$ with the tensor product $\widehat A\ot \widehat B$. Indeed, given $\omega=\psi_\#(\,\cdot\,a'b')$, with $a'\in A$ and $b'\in B$, we have
$$\omega(ba)=\psi_\#(baa'b')=\psi_A(aa')\psi_B(bb')$$
for all $a\in A$ and $b\in B$. This is of course much simpler but unfortunately, it is not what we need (see the proof of Theorem 3.7 below). The difference between the two identifications lies in the use of the adjoint of the map $R$, defined from $B\ot A$ to $A\ot B$ by
$$R(b\ot a)=\sum_{(a)}a_{(1)}\ot (b\tl a_{(2)}).$$
As explained earlier, it is far from obvious that this adjoint will map $\widehat A\ot \widehat B$ into $\widehat B\ot \widehat A$. We will investigate this result in the next item.
\nl
\it The right coaction of $C$ on $D$ \rm
\nl 
Recall that $C=\widehat A$ and $D=\widehat B$ and that formally, the right coaction $\Gamma$ is defined by the equation
$\langle b\ot a,\Gamma(d)\rangle=\langle b\tl a, d\rangle$ 
when $a\in A$, $b\in B$ and $d\in D$. This formula however only defines $\Gamma(d)\in (B\ot A)'$, the linear dual space of $B\ot A$. Similarly, we can use the equation 
$\langle b\ot a,(1\ot c)\Gamma(d)\rangle=\langle R(b \ot a), c\ot d\rangle$ 
to define a map $c\ot d\mapsto (1\ot c)\Gamma(d)$ from $C\ot D$ to $(B\ot A)'$. 
\snl
We will begin with showing, among other things, that this last equation does define a map from $C\ot D$ to $D\ot C$. 

\inspr{3.2} Proposition \rm
The map $R:B\ot A \to A\ot B$, defined before by 
$$R(b\ot a)=\sum_{(a)}a_{(1)} \ot (b\tl a_{(2)}) \tag"(3.1)"$$ 
has an adjoint $T:C \ot D \to D\ot C$, defined by
$$\langle b \ot a, T(c\ot d) \rangle=\langle R(b\ot a),c\ot d\rangle$$ 
for all $a\in A$, $b\in B$ $c\in C$ and $d\in D$. The map $T$ is again bijective.

\snl\bf Proof\rm:
First, we show that the above equation does define $T(c\ot d)$ in $D\ot C$. To do this, take $a'\in A$ and $b'\in B$ and let $c=\psi_A(a'\,\cdot\,)$ and $d=\psi_B(\,\cdot\, b')$. Then we find for all $a\in A$ and $b\in B$:
$$\align \langle R(b\ot a),c\ot d \rangle
	&=\sum_{(a)} \psi_A(a'a_{(1)})\psi_B((b\tl a_{(2)})b')\\
	&=\sum_{(a)} \psi_\#(a'a_{(1)}(b\tl a_{(2)})b')\\
	&=\psi_\#(a'bab')\\
	&=\psi_\#(bab'\sigma'_\#(a')).
\endalign$$
Because $b'\sigma'_\#(a')$ is in $AB$, we can write it as $\sum_i p_iq_i$ with $p_i\in A$ and $q_i\in B$, and we get
$$ \langle R(b\ot a),c\ot d\rangle = \sum_i\psi_A(ap_i)\psi_B(bq_i)$$
(where we have used the second formula in Proposition 2.2 in the previous section). So we can define $T(c\ot d)=\sum_i \psi_B(\,\cdot\,q_i)\ot \psi_A(\,\cdot\, p_i)$ in $\widehat B\ot \widehat A$. 
\snl
It is clear from the definition and the fact that $R$ is surjective, that $T$ is injective. It is also possible to show that $T$ is surjective. To prove this, one has to go once more through the argument and see that we can obtain all elements in $D\ot C$ in this way. One may need to use that $b'\sigma'_\#(a')=\sigma'_\#({\sigma'_B}^{-1}(b)a')$. 
\hfill $\square$
\einspr

In a completely similar way, one can show that also the map 
$$R^{\text{op}}:b\ot a\mapsto \sum_{(a)}a_{(2)}\ot (b\tl a_{(1)})\tag"(3.2)"$$
for $a\in A$ and $b\in B$ has a bijective adjoint $T^{\text{op}}:C\ot D \to D\ot C$. Here we have to start with $c=\psi_A(S_A(\,\cdot\,)a')$ and $d=\psi_B(\,\cdot\,b')$ where $a'\in A$ and $b'\in B$. 
\nl
Now, we are ready to obtain the first important result in this section.

\inspr{3.3} Theorem \rm
There exists a right coaction $\Gamma:D\to M(D\ot C)$, defined by 
$$\langle b\ot a,\Gamma(d)\rangle = \langle b\tl a, d\rangle $$
where $a\in A$, $b\in B$ and $d\in D$, making $D$ into a right $C$-comodule coalgebra.

\snl\bf Proof\rm:
As mentioned before, we can define $\Gamma(d)\in (B\ot A)'$ for $d\in D$ using this formula. We see however easily that then 
$$\langle b \ot a, (1\ot c)\Gamma(d) \rangle=\langle R(b\ot a),c\ot d\rangle$$ 
for all $a\in A$, $b\in B$ $c\in C$ and $d\in D$. And from Proposition 3.2, it follows that actually $(1\ot c)\Gamma(d)\in D\ot C$. Similarly, we find that $\Gamma(d)(1\ot c)\in D\ot C$ for all $c\in C$ and $d\in D$. Therefore $\Gamma(d)\in M(D\ot C)$ for all $d$. The injectivity of $\Gamma$ follows from the injectivity of $T$.
\snl
Now, we know from (a remark following) Proposition 1.2 in [De-VD-W] that $R$ satisfies the equations
$$\align R(\iota_B\ot m_A)
      &=(m_A\ot \iota_B)(\iota_A\ot R)(R\ot \iota_A) \qquad \text{on} \qquad B\ot A\ot A \\
   R(m_B\ot\iota_A) 
      &=(\iota_A\ot m_B)(R\ot\iota_B)(\iota_B\ot R) \qquad \text{on} \qquad B\ot B\ot A.
\endalign$$
Therefore, the adjoint $T$ will satisfy the adjoint equations
$$\align (\iota_D\ot\Delta_C)T
    	&=(T\ot \iota_C)(\iota_C\ot T)(\Delta_C\ot\iota_D) \qquad \text{on}\qquad C\ot D \\
         (\Delta_D\ot \iota_C)T
        &=(\iota_D\ot T)(T\ot \iota_D)(\iota_C\ot \Delta_D) \qquad \text{on}\qquad C\ot D
\endalign$$ 
needed for $D$ to be a  right $C$-comodule coalgebra (see Definition 1.13 in [De-VD-W]).

\hfill$\square$\einspr

We know that there is a covering difficulty with the interpretation of the last equations in the argument above. This has been discussed in [De-VD-W]. The problem can be overcome if we use the inverses of the maps $T\ot \iota_C$ and $\iota_D\ot T$ as is done in Remark 1.10 in [De-VD-W].
\snl
We now make some remarks about notations (in accordance with the notations of [De-VD-W]).

\inspr{3.4} Remark \rm
i) In general, we sometimes use the same symbols for objects related with the original pair $(A,B)$ and the corresponding ones for the dual pair $(C,D)$. We take into account however that these matched pairs are of a different kind.
\snl
ii) There will be no confusion with the action as we use $\tl$ for the right action of $A$ on $B$ and we will use $\tr$ for the left action of $D$ on $C$ (to be defined later). We use however $\Gamma$ for the coactions in the two cases. For the first pair $(A,B)$ we have $\Gamma: A\to M(B\ot A)$ whereas for the second pair $(C,D)$, we have $\Gamma:D\to M(D\ot C)$. See Theorem 3.3 above.
\snl
iii) We have the twist maps $R$ and $R^{\text{op}}$. For the original pair we have recalled the formulas already (cf.\ the formulas (3.1) and (3.2). They are maps from $B\ot A$ to $A\ot B$. For the new pair we will have maps from $D\ot C$ to $C\ot D$, given by 
$$\align R(d\ot c) &=\sum_{(d)} (d_{(1)}\tr c) \ot d_{(2)} \tag"(3.3)"\\
         R^{\text{op}}(d\ot c) &=\sum_{(d)} (d_{(2)}\tr c) \ot d_{(1)}.\tag"(3.4)"
\endalign$$
iv) We finally have a similar situation for the cotwist maps $T$ and $T^{\text{op}}$. For the original pair, these are maps from $A\ot B$ to $B\ot A$. Recall that they are given by
$$T(a\ot b)=\Gamma(a)(b\ot 1) \qquad\quad \text{and} \qquad\quad T^{\text{op}}(a\ot b)=(b\ot 1)\Gamma(a).\tag"(3.5)"$$
For the dual pair, we have maps from $C\ot D$ to $D\ot C$, given by
$$T(c\ot d)=(1\ot c)\Gamma(d) \qquad\quad \text{and} \qquad\quad T^{\text{op}}(c\ot c)=\Gamma(d)(1\ot c).\tag"(3.6)"$$
\hfill$\square$\einspr

Observe that in Proposition 3.2, we  {\it first define} $T$ on $C\ot D$ as the adjoint of the map $R$ on $B\ot A$. Similarly we first define $T^{\text{op}}$ on $C\ot D$ as the adjoint of the map $R^{\text{op}}$ on $B\ot A$. Only after Theorem 3.3, we can argue that these maps satisfy the formulas in (3.6) as expected. We will get a similar situation with the maps $R$ and $R^{\text{op}}$ on $D\ot C$. We will first define them as adjoints, use these to get the coaction and then argue that they are precisely given as in (3.3) and (3.4) - see the next item.
\snl
Using the same symbols for different maps should not be confusing here as we will be very systematic in the choice of the letters for elements in $A$, $B$, $C$  and $D$.
\nl
Before we continue, let us {\it relate} the existence of the coaction $\Gamma$ with a property found in Section 2. Indeed, with the terminology used in the previous section, we showed, (see a remark after) Lemma 2.7, that the map $b\mapsto \langle b\tl a, d \rangle$ will be again in $D$ for all $a\in A$ and $d\in D$. This means that we define a left action of $A$ on $D$ by $\langle b \tl a, d \rangle=\langle b,a\tr d \rangle$. Also this action is unital. What we have done here is going one step further in the duality and also dualize with respect to the variable $a\in A$.
\nl
\it The left action of $D$ on $C$ \rm
\nl
Now, we proceed along the same lines as before to get the left action of $D$ on $C$ that is dual to the left coaction of $B$ on $A$. We know already that formally 
$$\langle a,d\tr c\rangle=\langle \Gamma(a),d\ot c\rangle$$ 
for $a\in A$, $c\in C$ and $d\in D$. This defines $d\tr c$ in $A'$. Similarly, the equation  
$$\langle a\ot b ,R(d\ot c)\rangle = \langle \Gamma(a)(b\ot 1),d\ot c\rangle$$
for $a\in A$, $b\in B$, $c\in C$ and $d\in D$ will define $R(d\ot c)\in (A\ot B)'$. 
\snl
In this item, we show that these elements belong to the right spaces.
\snl
In the next proposition, we find the adjoint $R^{\text{op}}$ of the map $T^{\text{op}}$. It would of course be more natural to consider $T$ in stead of $T^{\text{op}}$ but it turns out that this is more complicated. We will say a bit more after the proof of this property.

\inspr{3.5} Proposition \rm 
The map $T^{\text{op}}:A\ot B\to B\ot A$, defined by $T^{\text{op}}(a\ot b)=(b\ot 1)\Gamma(a)$ for $a\in A$ and $b\in B$ has an adjoint $R^{\text{op}}:D\ot C\to C\ot D$, defined by 
$$\langle a\ot b ,R^{\text{op}}(d\ot c)\rangle = \langle (b\ot 1)\Gamma(a),d\ot c\rangle$$
for $a\in A$, $b\in B$, $c\in C$ and $d\in D$. This map $R^{\text{op}}$ is bijective.

\snl\bf Proof\rm:
Take $a'\in A$, $b'\in B$ and put $c=\psi_A(\,\cdot\,a')$ and $d=\psi_B(\,\cdot\,b')$. Also choose $a''\in A$ and $b''\in B$ so that $\psi_A(a'')=1$ and $\psi_B(b'')=1$.
For all $a\in A$ and $b\in B$ we have that
$$\align \langle(b\ot 1)\Gamma(a),d\ot c\rangle
	&=(\psi_B\ot\psi_A)((b\ot 1)\Gamma(a)(b'\ot a'))\\
	&=(\psi_\#\ot\psi_\#)((b\ot 1)\Gamma(a)(b'a''\ot a'b''))\\
	&=\sum_{(a)}(\psi_\#\ot\psi_\#)((bS(a_{(1)})\ot 1)\Delta_\#(a_{(2)})(b'a''\ot a'b'')).
\endalign$$
Now write $b'a''\ot a'b''=\sum_i \Delta_\#(x_i)(y_i\ot 1)$ with $x_i, y_i\in AB$ and assume that $x_i=p_iq_i$ with $p_i\in A$ and $q_i\in B$. Then we find
$$\align \langle(b\ot 1)\Gamma(a),d\ot c\rangle
	&=\sum_{i,(a)}\psi_\#(bS_A(a_{(1)})\delta_\#^{-1}y_i)\psi_\#(a_{(2)}x_i) \\
 	&=\sum_{i,(a)}\psi_\#(bS_A(a_{(1)})\delta_\#^{-1}y_i)\psi_A(a_{(2)}p_i)\psi_B(q_i) \\
	&=\sum_{i,(p_i)}\psi_\#(bp_{i(1)}\delta_A\delta_\#^{-1}y_i)\psi_A(ap_{i(2)})\psi_B(q_i).
\endalign$$
From this last expression, we see that we can define $R^{\text{op}}(d\ot c)$ in $C\ot D$ satisfying
$$\langle a\ot b ,R^{\text{op}}(d\ot c)\rangle = \langle (b\ot 1)\Gamma(a),d\ot c\rangle$$
for all $a,b$.
\snl
The injectivity of this map follows from the surjectivity of the map $T^{\text{op}}$. And to prove the surjectivity, we again have to go carefully through the proof and observe that every element in $C\ot D$ can be obtained. 	
\hfill $\square$ 
\einspr

As mentioned already, it would have been more natural to look at the map $T$ itself, defined from $A\ot B$ to $B\ot A$  by $T(a\ot b)=\Gamma(a)(b\ot 1)$. Unfortunately, a similar argument as above for the adjoint of this map, does not seem to work. It can be shown however that also the adjoint $R:D\ot C\to C\ot D$ of $T$ exists and is again bijective. One possibility is to use  the equality $T\circ R=T^{\text{op}}\circ R^{\text{op}}$, the fact that both $T^{\text{op}}$ and $R^{\text{op}}$ are adjointable so that also $T\circ R$ is adjointable and finally, that $R$ has a bijective adjoint. This will yield the adjoint of $T$ and it will still be bijective. 
\snl
In fact, we do not really need this argument as it will also follow from the second main result we prove now. 

\inspr{3.6} Theorem \rm 
There exists a left action $\tr$ of $D$ on $C$, defined by 
$$\langle a,d\tr c\rangle=\langle \Gamma(a),d\ot c\rangle$$ 
with $a\in A$, $c\in C$ and $d\in D$, making $C$ into a left $D$-module algebra.

\snl \bf Proof\rm:
We have mentioned already that the above formula will define $d\tr c$ in $A'$ for all $c\in C$ and $d\in D$. Now, because we have 
$$\langle(b\ot 1)\Gamma(a),d\ot c \rangle = \sum_{(d)}\langle a\ot b, (d_{(2)}\tr c) \ot d_{(1)} \rangle,$$
it follows from the previous proposition that 
$$\sum_{(d)}(d_{(2)}\tr c) \ot d_{(1)} \in C\ot D $$
for all $c\in C$ and $d\in D$. By applying the counit $\varepsilon_D$ on the last factor, we get $d\tr c\in C$ for all $c\in C$ and $d\in D$.
\snl
Next, if we dualize the appropriate equations for $T^{\text{op}}$ to the necessary equations for its adjoint $R^{\text{op}}$, just as we did in the proof of Theorem 3.3, we see that it will follow that $C$ is a left $D$-module algebra because $A$ is a left $B$-comodule coalgebra. This completes the proof.
\hfill $\square$ 
\einspr

Again, now that we have obtained the left action of $D$ on $C$, we can show that the adjoints $R$ and $R^{\text{op}}$ of the maps $T$ and $T^{\text{op}}$ indeed are given by the formulas (3.3) and (3.4) as claimed in Remark 3.4.iii.
\snl
We are now ready for the main subsection of this part.
\nl
\it The dual of $AB$ is $CD$ \rm
\nl
So, given the pair $(A,B)$ of algebraic quantum groups with the right action $\tl$ of $A$ on $B$ and the left coaction $\Gamma$ of $B$ on $A$, we can associate a right coaction $\Gamma$ of $C$ on $D$ and a left action $\tr$ of $D$ on $C$ by duality. This is proven in Theorem 3.3 and Theorem 3.6 above.
\snl
We can consider the smash product $CD$ and the smash coproduct $\Delta_\#$ on $CD$. The smash product is defined using the map $R$ as given in (3.3) on $D\ot C$. The algebra $CD$ is the algebra generated by $C$ and $D$ with the commutation rules 
$$dc=\sum_{(d)}(d_{(1)}\tr c)d_{(2)}$$
for $c\in C$ and $d\in D$. The coproduct $\Delta_\#$ on $CD$ is given by the formula
$$\Delta_\#(cd)=\sum_{(c),(d)}(c_{(1)}\ot c_{(2)})\Gamma(d_{(1)})(1\ot d_{(2)})$$
again for $c\in C$ and $d\in D$. Remember that the formulas for the pair $(C,D)$ can be found by either using the 'conversion' procedure or by 'duality' (as explained in Section 1 [De-VD-W]).
\snl
Then, we have the following {\it main result} of this section.

\inspr{3.7} Theorem \rm
The left action of $D$ on $C$ and the right coaction of $C$ on $D$ make $(C,D)$ into a matched pair (of the second type) and the bicrossproduct $(CD,\Delta_\#)$ is identified with the dual $(AB)\,\widehat{}$ of $AB$ as algebraic quantum groups by means of the natural pairing between $AB$ and $CD$ given by $\langle ab,cd \rangle=\langle a,c \rangle \langle b,d\rangle$.

\snl \bf Proof\rm:
We first have to verify that $(C,D)$ is a matched pair (of the second type). This means that the operator $P$, defined on $D\ot C$ as $T^{\text{op}}\circ R^{\text{op}}$, has to satisfy the equations $P=T \circ R$ and
$$\align 
    P(m_D\ot \iota_C) &= (m_D\ot \iota_C)P_{13}P_{23} \quad \text{on} \quad D\ot D\ot C \\
   (\iota_D\ot \Delta_C) P &= P_{12}P_{13}(\iota_D\ot\Delta_C) \quad \text{on} \quad D\ot C
\endalign$$
of Theorem 2.16 in [De-VD-W]. This follows by duality. Indeed, by assumption, these (or rather similar) equations are satisfied for the map $P$ on $B\ot A$ as $(A,B)$ is assumed to be a matched pair (as in Theorem 2.14 of [De-VD-W]).  Because both the operators $T\circ R$ and $T^{\text{op}} \circ R^{\text{op}}$, are  self-dual, the equations follow by duality. So, we do get that the smash coproduct $\Delta_\#$ makes the smash product $CD$ into an algebraic quantum group by Theorem 2.16 of [De-VD-W].
\snl
Now, we use the result of Proposition 3.1 of this section, giving the identification of the dual $(AB)\,\widehat{}\,$ with $CD$ as linear spaces. This is realized with the pairing 
$$\langle ab,cd \rangle=\langle a,c \rangle \langle b,d\rangle$$
where $a\in A$, $b\in B$, $c\in C$ and $d\in D$.
\snl
We only have to argue that the products and coproducts are dual to each other. But this is again quite obvious as the pairing is just the tensor product pairing of $A\ot B$ with $C\ot D$ and because e.g.\ the product on $A\ot B$ is given by the formula
$$(m_A\ot m_B)(\iota_A\ot R\ot \iota_B)$$
(see Proposition 1.3 in [De-VD-W]) and the coproduct on $C\ot D$ is given by 
$$(\iota_C\ot T\ot \iota_D)(\Delta_C\ot \Delta_D)$$
(see Section 2 in [De-VD-W]). Clearly, these two definitions are dual to each other. This completes the proof. 
\hfill $\square$\einspr  

Before we continue with a brief look at the $^*$-algebra case, we first come back to a result announced in Section 1. In Proposition 1.1 we have the homomorphisms $\rho$ and $\eta$. We have shown in Proposition 1.3 that these are equal if $B$ is assumed to be a right $A$-module {\it bi}-algebra. Now, we can use duality and the result of Proposition 1.5 to show that this is also true in general.

\inspr{3.8} Corollary \rm Let $(A,B)$ be as before in this section. Assume now that there are cointegrals in $B$. Let $h$ be a left cointegral and let $k$  a right cointegral. Then there is a (single) homomorphism  $\eta:A\to \Bbb C$ such that
$$h\tl a=
\eta(a)h \qquad\qquad \text{and} \qquad\qquad k\tl a = \eta(a)k$$
for all $a\in A$.
\snl\bf Proof\rm:
We know by applying the results of Proposition 1.4 and 1.5 to the dual case that there exists an element $y\in M(C)$ satisfying
$$(\varphi_D\ot\iota)\Gamma(d)=\varphi_D(d)y\qquad\qquad\text{and}\qquad\qquad (\psi_D\ot\iota)\Gamma(d)=\psi_D(d)y$$
for all $d\in D$ (where we have used the notations as earlier in this section). We have $\varphi_D(d)=\langle h,d\rangle$ and $\psi_D(d)=\langle k,d\rangle$ for $d\in D$ (provided the cointegrals are properly normalized). Then, pairing these equations with any element $a\in A$, we find
$$\align \langle h\tl a,d\rangle 
          &=\langle h \ot a, \Gamma(d)\rangle = \langle h,d\rangle \langle y ,a\rangle \\
         \langle k\tl a,d\rangle 
          &=\langle k \ot a, \Gamma(d)\rangle = \langle k,d \rangle \langle y,a\rangle 
\endalign$$
and we see that $\eta(a)=\langle y,a\rangle$ for all $a\in A$.
\hfill $\square$\einspr

 \it The $^*$-algebra case \rm
 \nl
We now look briefly at the case where $A$ and $B$ are $^*$-algebras. Then also the duals $C$ and $D$ are $^*$-algebras and the involutions satisfy
$$\langle a,c^* \rangle= \langle S_A(a)^*, c\rangle^-\qquad\quad \text{and} \qquad\quad
\langle b ,d^* \rangle=\langle S_B(b)^*,d \rangle^-$$
for all $a\in A, b\in B$ and $c\in C, d\in D$.
\snl
The relevant results in this case are given in the following proposition.

\inspr{3.9} Proposition \rm
Assume that the algebras $A$ and $B$ are $^*$-algebras and that the right action of $A$ on $B$ and the left coaction of $B$ on $A$ satisfy
$$(b\tl a)^*=b^*\tl S_A(a)^*
\qquad\qquad \text{and} \qquad\qquad
\Gamma(S_A(a)^*)=((\iota\ot S_A)\Gamma(a))^*\tag"(3.7)"$$
for all $a\in A$, $b\in B$ and $d\in D$ (as in Theorem 2.15 of [De-VD-W]). Then also the dual left action of $D$ on $C$ and the dual right coaction of $C$ on $D$ satisfy
 $$(d\tr c)^*=S(d)^*\tr c^*
\qquad\qquad \text{and} \qquad\qquad
\Gamma(S_D(d)^*)=((S_D\ot\iota)\Gamma(d))^*\tag"(3.8)"$$
for all $a\in A$, $c\in C$ and $d\in D$. Moreover, the involution on $CD$, obtained as the dual of $AB$ coincides with the involution on $CD$ given by $(cd)^*=d^*c^*$ for all $c\in C$ and $d\in D$.

\snl\bf Proof\rm:
The proof is rather straightforward. Take e.g.\ $a\in A$, $c\in C$ and $d\in D$. Then
$$\align \langle a, (d\tr c)^* \rangle &= \langle S_A(a)^*,d\tr c \rangle^- \\
&=\langle \Gamma(S_A(a)^*),d\ot c \rangle^- \\
&=\langle ((\iota\ot S_A)\Gamma(a))^*,d\ot c \rangle^- \\
&=\langle \Gamma(a),S_D(d)^*\ot c^* \rangle \\
&=\langle a, S_D(d)^*\tr c^* \rangle.
\endalign$$
This proves the first formula in (3.8). The proof of the second one is using completely the same arguments (in the other order and for the dual system).
\snl
Finally, to prove the last statement, take elements $a,b,c,d$. Using the definition of $S_\#$ as found in Theorem 3.1 of [De-VD-W], we get
$$\align \langle ab, (cd)^* \rangle &= \langle S_\#(ab)^*,cd \rangle^- \\
&=\sum_{(a)}\langle S_A(a_{(0)})^* S_B(a_{(-1)}b)^*  ,cd \rangle^- \\
&=\sum_{(a)}\langle S_A(a_{(0)})^*,c\rangle^- \langle S_B(a_{(-1)}b)^*  ,d \rangle^- \\
&=\sum_{(a)}\langle a_{(0)},c^*\rangle \langle a_{(-1)}b  ,d^* \rangle \\
&=\langle T(a\ot b),d^*\ot c^* \rangle \\
&=\langle a\ot b, R(d^*\ot c^* ) \rangle= \langle ab, d^*c^*\rangle.
\endalign$$ 
This completes the proof.
\hfill $\square$\einspr

If we look at the above proof, we see (again) why the conditions (3.7) and (3.8)  are natural. Moreover, from the last argument above, we also see where the formula for $S_\#$ on $AB$ comes from. Finally, in Proposition 2.12 we have shown that $AB$ has positive integrals if $A$ and $B$ have positive integrals. We know that in this case, also the duals have positive integrals.
\nl
We will not consider examples in this section but refer to the next special section on examples for illustrations of various results in this section.
\nl\nl

\bf 4. Examples and special cases\rm
\nl
In this section, we will try to illustrate as many of the results in this paper as possible, with various examples. At the same time, we look at special cases and consider examples within these cases.
\snl
It is most natural to start with the simplest example, coming from a matched pair of (discrete, but possibly infinite) groups, as discussed already in detail in the first paper [De-VD-W]. Unfortunately, this is too simple. Most, if not all of the objects discussed in this paper are trivial for that example. The reason is of course that left and right integrals coincide and that the algebras are either abelian or that the integrals are traces. This forces the modular elements, as well as the modular automorphisms, scaling constants, ... for the components to be trivial. Then it is expected that this will also be the case for the bicrossproducts. We will shortly verify this.
\nl
\it The example of a matched pair of groups \rm

\inspr{4.1} Example \rm 
Recall that we consider a group $G$ with two subgroups $H$ and $K$ so that $G=KH$ and $H\cap K=\{e\}$ where $e$ is the identity in the group. We  have a left action $\tr$ of $H$ on $K$ and a right action $\tl$ of $K$ on $H$ given by
$hk=(h\tr k)(h\tl k)$ whenever $h\in H$ and $k\in K$. The algebra $A$ is the group algebra $\Bbb C H$ and the algebra $B$ is the algebra $F(K)$ of complex functions with finite support on $K$. The right action of $A$ on $B$ is given by the formula $(f\tl h)(k)=f(h\tr k)$ when $f\in B$, $h\in H$ and $k\in K$. The left coaction $\Gamma$ of $B$ on $A$ is given by $\Gamma(h)(\delta_k\ot 1)=\delta_k \ot (h\tl k)$ when $h\in H$ and $k\in K$ and where $\delta_k$ is the function in $B$ that takes the value $1$ in the element $k$ and $0$ everywhere else. We consider the group $H$ as sitting in the group algebra $\Bbb C H$.
\snl
A left integral $\varphi_A$ is given by $\varphi_A(h)=0$ when $h\in H$ and $h\neq e$ and $\varphi_A(e)=1$. It is also right invariant (so that $\delta_A=1$) and it is a trace (so that the modular automorphisms are all trivial). A right integral $\psi_B$ on $B$ is given $\psi_B(f)=\sum_k f(k)$, it is also left invariant (so that also $\delta_B=1$) and as the algebra $B$ is abelian, again the modular automorphisms are trivial. Of course, the two scaling constants are $1$ because the square of the antipode is the identity map.
\snl
Because $h\tl k=e$ if and only if $h=e$, we find that
$$(\iota\ot\varphi_A)\Gamma(h)=\varphi_A(h)1$$
for all $h\in H$ and therefore, the element $y$, as defined in Proposition 1.4, is equal to $1$ in this example. This will give that also $\delta_\#=1$. Furthermore 
$$\psi_B(f\tl h)=\sum_k f(h\tr k)=\sum_k f(k)=\psi_B(f)$$
when $f\in B$ and $h\in H$. This means that the map $\gamma$, from Proposition 2.9, coincides with $a\mapsto \varepsilon_A(a)1$ and this implies that also the modular automorphisms $\sigma_\#$ and $\sigma'_\#$ are trivial. 
\snl
Because the multiplier $y$ is trivial, we have for the left and right integrals $\varphi_\#$ and $\psi_\#$ on the bicrossproduct simply
$$\varphi_\#(ab)=\varphi_A(a)\varphi_B(b) 
\qquad\quad \text{and}\qquad\quad \psi_\#(ab)=\psi_A(a)\psi_B(b)$$
whenever $a\in A$ and $b\in B$.
\snl
Let us finish this example with a brief look at duality. Of course the dual $C$ of $A$ is the function algebra $F(H)$ whereas the dual $D$ of $B$ is the group algebra $\Bbb C K$ of $K$. Moreover, the left action of $D$ on $C$, as obtained by duality in Theorem 3.6 is given, as expected, by $k\tr f=f(\,\cdot\,\tl k)$ and the right coaction of $C$ on $D$, obtained by duality in Theorem 3.3 is given by $\Gamma(k)(1\ot \delta_h)=(h\tr k)\ot \delta_h$. We recover the dual bicrossproduct as given e.g.\ in Example 2.17 of [De-VD-W]. 
\hfill$\square$
\einspr
 
We see that indeed this example is far too simple to illustrate our results. Nevertheless, we will come back to this case at the end of the section to illustrate some other phenomena (see Proposition 4.18).
\nl
\it The example with $(B,\Delta_B)=(A,\Delta_A^{\text{cop}})$ \rm
\nl
Let us now consider once more Example 1.10 that we have used before, throughout the various sections in the paper. 
\snl
Recall that we have a regular multiplier Hopf algebra $(A,\Delta_A)$ and that we have put $(B,\Delta_B)=(A,\Delta_A^{\text{cop}})$. The action is the adjoint action, defined by
$$b\tl a=\sum_{(a)}S(a_{(1)})ba_{(2)}$$
and the coaction is the adjoint coaction given by
$$\Gamma(a)=\sum_{(a)}S(a_{(1)})a_{(3)}\ot a_{(2)}.$$
Remember that when using the Sweedler notation here, we use it for the original coproduct $\Delta_A$. In what follows, the antipode $S$ is the antipode $S_A$ of $A$. Then $S_B=S^{-1}$.
\snl
Just as for ordinary Hopf algebras, we have the following result. The proof is essentially the same, but we include it for completeness. Again, as we did before, we will use again the notation $A\# B$ for the bicrossproduct. 

\inspr{4.2} Proposition \rm 
Take $(A,\Delta_A)$ and  $(B,\Delta_B)$ as above. Define $\theta:A\ot B \to A\# B$ by
$$\theta (a\ot b)=\sum_{(a)}a_{(1)} \# S(a_{(2)})b.$$
Then $\theta$ is an isomorphism from the tensor product of the multiplier Hopf algebras $A$ and $B$ with the bicrossproduct $(A\# B,\Delta_\#)$.

\snl\bf Proof\rm: Take $a,a'\in A$ and $b,b'\in B$. We have on the one hand
$$\align \theta((a\ot b)(a'\ot b'))=\theta(aa'\ot bb')
           &=\sum_{(a)}a_{(1)}a'_{(1)}\# S(a_{(2)}a'_{(2)})bb' \\
           &=\sum_{(a)}a_{(1)}a'_{(1)}\# S(a'_{(2)})S(a_{(2)})bb' 
\endalign$$
while on the other hand, we have
$$\align \theta(a\ot b)\theta(a'\ot b')
           &=\sum_{(a),(a')}(a_{(1)}\# S(a_{(2)})b)(a'_{(1)}\# S(a'_{(2)})b') \\
           &=\sum_{(a),(a')}a_{(1)}a'_{(1)}\# ((S(a_{(2)})b)\tl a'_{(2)})S(a'_{(3)})b'\\
           &=\sum_{(a),(a')}a_{(1)}a'_{(1)}\# S(a'_{(2)})S(a_{(2)})ba'_{(3)}S(a'_{(4)})b'\\
           &=\sum_{(a),(a')}a_{(1)}a'_{(1)}\# S(a'_{(2)})S(a_{(2)})bb'.
\endalign$$
These two expressions are the same, proving that $\theta$ is an algebra map. Because it is clearly bijective, we see that we have an isomorphism of algebras.
\snl
Next, we show that this isomorphism converts the coproduct $\Delta_\#$ on $A\# B$ to the tensor coproduct on $A\ot B$. First one can verify that 
$$\Delta_\#(a\# b)=\sum_{(a)(b)}(a_{(1)}\# S(a_{(2)})a_{(4)}b_{(2)})\ot (a_{(3)}\# b_{(1)})$$
whenever $a\in A$ and $b\in B$. Remember that $\Delta_B(b)=\sum_{(b)} b_{(2)}\ot b_{(1)}$ because we use the Sweedler notation for $\Delta_A$. Then we find
$$\align \Delta_\#(\theta(a\ot b)) 
          &=\sum_{(a)}\Delta_\#(a_{(1)}\# S(a_{(2)})b) \\
          &=\sum_{(a),(b)}(a_{(1)}\# S(a_{(2)})a_{(4)}S(a_{(5)})b_{(2)})\ot (a_{(3)}\# S(a_{(6)})b_{(1)})\\
          &=\sum_{(a),(b)} (a_{(1)}\# S(a_{(2)})b_{(2)})\ot (a_{(3)}\# S(a_{(4)})b_{(1)})\\
          &=\sum_{(a),(b)} \theta(a_{(1)} \ot b_{(2)})\ot \theta(a_{(2)}\ot b_{(1)})
\endalign$$
and we see that $\theta$ is also a coalgebra map.
\hfill$\square$
\einspr

Because of this result, it is {\it not really necessary} to verify that the action and the coaction are matched. It follows because we know from the above proposition that $\Delta_\#$ is a homomorphism on $A\# B$. We refer to a remark made in Section 1, when treating Example 1.10.
\snl
It also follows that the bicrossproduct {\it gives nothing really new}. On the other hand, precisely because of this result, it is easy and interesting to check the obtained results for this example. This is what we will do next.

\inspr{4.3} Example \rm
We continue the investigation we started in Example 4.2.
\snl
Let us e.g.\ verify the formula for $S_\#$. In Theorem 3.1 of [De-VD-W] we find that, in the general case  
$$S_\#(ab)=\sum_{(a)}S_B(a_{(-1)}b)S_A(a_{(0)})$$
when writing the smash product as $AB$. However, in the case of our example, we can not use this notation and we have to be more careful. In this case we consider first
$$\sum_{(a)}a_{(-1)}b\ot a_{(0)}=\sum_{(a)}S(a_{(1)})a_{(3)}b\ot a_{(2)}$$
and so 
$$\align \sum_{(a)}S_B(a_{(-1)}b)\ot S_A(a_{(0)})&=\sum_{(a)}S^{-1}(S(a_{(1)})a_{(3)}b)\ot S(a_{(2)})\\
          &=\sum_{(a)}S^{-1}(b)S^{-1}(a_{(3)})a_{(1)}\ot S(a_{(2)}).
\endalign$$
Next we apply the twist map $R$ and get
$$\align
R(\sum_{(a)}S_B(a_{(-1)}b)\ot S_A(a_{(0)}))
&=\sum_{(a)}S(a_{(3)})\ot (S^{-1}(b)S^{-1}(a_{(4)})a_{(1)})\tl S(a_{(2)})\\
&=\sum_{(a)}S(a_{(4)})\ot S^2(a_{(3)})S^{-1}(b)S^{-1}(a_{(5)})a_{(1)}S(a_{(2)})\\
&=\sum_{(a)}S(a_{(2)})\ot S^2(a_{(1)})S^{-1}(b)S^{-1}(a_{(3)}).
\endalign$$
It follows that 
$$S_\#(a\# b)=\sum_{(a)}S(a_{(2)})\# S^2(a_{(1)})S^{-1}(b)S^{-1}(a_{(3)}).$$
Then 
$$\align S_\#(\theta(a\ot b))&=\sum_{(a)}S_\#(a_{(1)}\#S(a_{(2)})b)\\
&=\sum_{(a)}S(a_{(2)})\#S^2(a_{(1)})S^{-1}(S(a_{(4)})b)S^{-1}(a_{(3)})\\
&=\sum_{(a)}S(a_{(2)})\#S^2(a_{(1)})S^{-1}(b)
\endalign$$
while also
$$\theta(S(a)\ot S^{-1}(b))=\sum_{(a)}S(a_{(2)})\#S^2(a_{(1)})S^{-1}(b).$$
This verifies the formula for $S_\#$ in this example. 
\hfill$\square$\einspr

Next, we consider the case where the original multiplier Hopf algebra has integrals. We know that then also the bicrossproduct has integrals and we can consider the various associated data. This is what we do in the next example.

\inspr{4.4} Example \rm
i) When we consider the formulas obtained for $\psi_\#$ and $\varphi_\#$ in Section 2, we find for this case
$$\align \psi_\#(\theta(a\ot b)) 
        &=\sum_{(a)}\psi_\#(a_{(1)}\# S(a_{(2)})b) \\
        &=\sum_{(a)}\psi_A(a_{(1)})\psi_B(S(a_{(2)})b) \\
        &=\psi_A(a)\psi_B(b)
\endalign$$
for all $a\in A$ and $b\in B$. Similarly
$$\align \varphi_\#(\theta(a\ot b)) 
        &=\sum_{(a)}\varphi_\#(a_{(1)}\# S(a_{(2)})b) \\
        &=\sum_{(a)}\varphi_A(a_{(1)})\varphi_B(yS(a_{(2)})b) \\
        &=\varphi_A(a)\varphi_B(\delta_AS(\delta_A)b)\\
        &=\varphi_A(a)\varphi_B(b)
\endalign$$
for all $a\in A$ and $b\in B$. We have used that $y=\delta_A$ as we have seen in Example 1.10.
\snl
ii) Next consider $\delta_\#$. From  the result in Proposition 2.3 (or rather from the remark following it), we find $\delta_\#=\delta_A\# \delta_A^{-1}\delta_B$. This is the same as $\theta(\delta_A\ot \delta_B)$ as expected.
\snl
iii) Finally, look at the modular automorphism $\sigma'_\#$. We have seen in Example 2.13 that
$$\sigma'_\#(a\#b)=\sum_{(a)}\sigma'_A(a_{(1)})\# S^{-1}(a_{(2)})\sigma_A(a_{(3)})\sigma'_B(b)$$
and so 
$$\align \sigma'_\#(\theta(a\ot b))
            &=\sum_{(a)}\sigma'_\#(a_{(1)}\# S(a_{(2)})b)\\
            &=\sum_{(a)}\sigma'_A(a_{(1)}) \# S^{-1}(a_{(2)})\sigma_A(a_{(3)})\sigma_A(S(a_{(4)}))\sigma_B'(b)\\
            &=\sum_{(a)}\sigma'_A(a_{(1)}) \# S^{-1}(a_{(2)})\sigma_B'(b)\\
            &=\sum_{(\sigma'_A(a))}(\sigma'_A(a))_{(1)}\# (S(\sigma'_A(a))_{(2)})\sigma_B'(b)\\ 
            &=\theta(\sigma'_A(a)\ot \sigma'_B(b)).
\endalign$$
In the above calculation, we have used that $\sigma'_B=\sigma_A$ because the coproduct $\Delta_B$ on $B$ is $\Delta_A^{\text{cop}}$ so that $\psi_B=\varphi_A$. We also have used that $\Delta_A\circ \sigma_A'=(\sigma_A'\ot S^{-2})\circ \Delta_A$. 
\hfill$\square$
\einspr

We see that, as it should, all the data that we found for this example convert to the right things under the automorphism $\theta$. Observe that for the scaling constant, we get $\tau_\#=\tau_A\tau_B$ and this is also the scaling constant for the tensor product $A\ot B$.
\nl
\it The duality for this example \rm
\nl
Consider again the situation discussed in the previous item. So we again start with a regular multiplier Hopf algebra $(A,\Delta_A)$ and with $(B,\Delta_B)=(A,\Delta_A^{\text{cop}})$. We assume that they have integrals and that we can consider the duals. We use $(C,\Delta_C)$ for the dual of $(A,\Delta_A)$ and $(D,\Delta_D)$ for the dual of $(B,\Delta_B)$ as in Section 3 on duality. Then we get that $(D,\Delta_D)=(C^{\text{op}},\Delta_C)$ because $(B,\Delta_B)=(A,\Delta_A^{\text{cop}})$. In this case, there is no confusion in using the coproduct and the Sweedler notation. For the product however, we {\it systematically use the product in} $C$. When we write the antipode $S$, we always mean the antipode $S_C$.
\snl
Then we get the following result.

\inspr{4.5} Proposition \rm
The adjoint of the right action of $A$ on $B$ (as in Theorem 3.3) gives the right coaction $\Gamma$ of $C$ on $D$ defined by
$$\Gamma(d)=\sum_{(d)} d_{(2)}\ot S(d_{(1)})d_{(3)}$$
for $d\in D$. Similarly, the adjoint of the left coaction of $B$ on $A$ (as in Theorem 3.6) gives the left action of $D$ on $C$, defined as
$$d\tr c=\sum_{(d)}S(d_{(1)})cd_{(2)}$$
for $c\in C$ and $d\in D$. Therefore, an application of Theorem 3.7 gives that the dual of $A\# B$ is the bicrossproduct $C\# D$, obtained from the matched pair $(C,D)$ as in Theorem 2.16 in [De-VD-W].
\hfill$\square$
\einspr

The proof is rather straightforward.
\snl
For finite-dimensional Hopf algebras, the bicrossproduct construction $C\# D$ in Proposition 4.5 is found as Example 6.2.2 of [M] where it is called the {\it mirror} construction, denoted by $M(C)$, of the Hopf algebra $C$. 
\snl
It is interesting to see what happens with the isomorphism $\theta$ in Proposition 4.2 in this duality. We get the following expected result. We use the notations as in the foregoing part.

\inspr{4.6} Proposition \rm 
Define $\eta:C\ot D\to C\# D$ by 
$$\eta(c\ot d)=\sum_{(d)}cd_{(1)}\# d_{(2)}$$
whenever $c\in C$  and $d\in D$. Then $\eta$ is an isomorphism from the tensor product of the multiplier Hopf algebras $C$ and $D$ with the bicrossproduct $C\# D$. Moreover we have
$$\langle \theta(a\ot b),\eta(c\ot d)\rangle=
\langle a\ot b,c\ot d\rangle=\langle a,c\rangle\langle  b ,d\rangle \tag"(4.1)")$$
for all $a,b,c,d$ in $A,B,C,D$ respectively.

\snl \bf Proof\rm:  The fact that $\eta$ is an isomorphism as stated in the proposition can be proven in a completely similar way as was done for the isomorphism $\theta$ in the proof of Proposition 4.2. In fact, it also follows from the formula (4.1). So, it is sufficient to prove that this formula holds. It is easy to verify:
$$\align
\sum_{(a),(d)}\langle a_{(1)} \#  S(a_{(2)})b, c d_{(1)}\# d_{(2)}\rangle
&=\sum_{(a),(d)}\langle  a_{(1)},c\rangle \langle  a_{(2)},d_{(1)}\rangle \langle  S(a_{(3)})b ,d_{(2)}\rangle \\
&=\sum_{(a),(d)}\langle  a_{(1)},c\rangle \langle a_{(2)} S(a_{(3)})b ,d\rangle \\
&=\langle  a,c\rangle \langle b ,d\rangle.
\endalign$$
This completes the proof as the second equality in (4.1) is true by definition.
\hfill$\square$
\einspr

\it The special case where the right action of $A$ on $B$ is trivial \rm
\nl
Consider a matched pair $(A,B)$ of regular multiplier Hopf algebras. We now moreover assume that {\it the right action $\tl$ of $A$ on $B$ is trivial}, that is that $b\tl a=\varepsilon_A(a)b$ for all $a\in A$ and $b\in B$. This implies $A$ and $B$ commute within the smash product $AB$ (because the twist map is nothing else but the flip map on $B\ot A$). In fact $AB$ is isomorphic with $A\ot B$ as an algebra. 
\snl
It is an easy consequence of the axioms that the coaction must satisfy $\Gamma(aa')=\Gamma(a)\Gamma(a')$ for all $a,a'\in A$. This can be seen e.g.\ from the fact that the coproduct $\Delta_\#$ is a homomorphism on $A$ and that $A$ and $B$ commute. See also a remark, following Theorem 2.16 in [De-VD-W]. So, we need to have that $A$ is a left $B$-comodule {\it bi}-algebra.
\snl
Moreover, as $R=R^{\text{op}}$ in this case the basic assumption that $T\circ R=T^{\text{op}}\circ R^{\text{op}}$ gives $T=T^{\text{op}}$ so that $\Gamma(a)(b\ot 1)=(b\ot 1)\Gamma(a)$ for all $a$ and $b$. This condition is of course fulfilled if also $\Gamma$ is trivial (i.e.\ when $\Gamma(a)=1\ot a$ for all $a$). It is also satisfied if $B$ is abelian. The first case is not interesting as then the bicrossproduct is nothing else but the tensor product of the two multiplier Hopf algebras. 
\snl
It is not difficult to get examples with $B$ abelian. The case where $B$ is abelian is in fact already considered in [De2]. Such examples are typically constructed from a matched pair of groups $(H,K)$ where the left action of $H$ on $K$ is trivial so that on the other hand, $K$ acts from the right on $H$ by means of isomorphisms of $H$. Then $G$ is a semi-direct product.
\snl
Let us now show that there are other cases where $T=T^{\text{op}}$.

\inspr{4.7} Example \rm 
Consider regular multiplier Hopf algebras $A,B_0,B_1$. Assume that $B_1$ is abelian and that $\Gamma_1: A\to M(B_1\ot A)$ is a left coaction of $B_1$ on $A$ making $A$ into a left $B_1$-comodule {\it bi}-algebra. Define $B=B_0\ot B_1$ with the tensor coproduct and $\Gamma:A\to M(B\ot A)$ by
$$\Gamma(a)=1\ot\Gamma_1(a)$$ for $a$ in $A$ where $1$ is the identity in $M(B_0)$ and where we view $M(B_0)\ot M(B_1\ot A)$ as sitting in $M(B\ot A)$. Then, it is not hard to see that still $A$ is a left $B$-comodule bi-algebra. Moreover we have that $b\ot 1$ and $\Gamma(a)$ will commute in $M(B\ot A)$ for all $a\in A$ and $b\in B$. Then, with the trivial right action of $A$ on $B$, we get a matched pair of multiplier Hopf algebras.
\hfill$\square$
\einspr

Of course, also this example is not very special, but at least, it shows that we do not need a trivial coaction of $B$ on $A$ or an abelian algebra $B$ in order to get a matched pair with a trivial action of $A$ on $B$. We are convinced that it should not be too difficult to get similar examples where $B$ is a twisted tensor product of two factors.
\nl
Next, we consider {\it this case with integrals}.

\inspr{4.8} Proposition \rm
Assume that $(A,B)$ is a matched pair of multiplier Hopf algebras with integrals and that the right action of $A$ on $B$ is trivial as before. Then, the multiplier $y$, as obtained in Proposition 1.4, is central in $M(B)$. 
\hfill$\square$
\einspr

This is an immediate consequence of Proposition 1.7 where it is shown that $y^{-1}b y = b\tl \delta_A$.
\snl
It is not completely obvious that there are cases where $y\neq 1$. We will give such an example below (see Example 4.12). First we consider the modular automorphisms.
\snl
We know from Proposition 2.6 that the modular automorphism $\sigma_\#$ and $\sigma'_\#$ on $B$ always coincide with the original modular automorphisms $\sigma_B$ and $\sigma'_B$. In this special case, this is also true for these modular automorphisms on $A$ as will follow from the next result.

\inspr{4.9} Proposition \rm Let $(A,B)$ be as in the previous proposition. Then we see that the linear map $\gamma$, defined from $A$ to $M(B)$ in Proposition 2.9, is given by $\gamma(a)=\varepsilon_A(a)1$. It follows that also $\sigma_\#$ and $\sigma'_\#$ on $A$ always coincide with the original modular automorphisms $\sigma_A$ and $\sigma'_A$.

\snl\bf Proof\rm:
It follows from the formula in Proposition 2.10 that the result is true for $\sigma'_\#$. Now, because 
$$\sigma_\#(a)=\delta_\#^{-1}\sigma'_\#(a)\delta_\#
=\delta_B^{-1} y \delta_A^{-1}\sigma'_\#(a)\delta_A y^{-1} \delta_B
=\delta_B^{-1} y \sigma_A(a)y^{-1} \delta_B$$
and because $A$ and $B$ commute, we see that also $\sigma_\#(a)=\sigma_A(a)$.
\hfill$\square$
\einspr

We now look for more concrete examples for this case. The following result is standard.

\iinspr{4.10} Proposition \rm
Let $(A,\Delta_A)$ be any multiplier Hopf algebra. Let $G$ be a group and consider an action $p\mapsto \alpha_p$ of $G$ by automorphisms of $(A,\Delta_A)$. Let $B$ be the algebra $K(G)$ of complex functions on $G$ with finite support and equip $B$ with the natural coproduct, coming from the group multiplication. Define $\Gamma: A\to M(B\ot A)$ by
$$\Gamma(a)(p)=\alpha_{p^{-1}}(a)$$
where $a\in A$ and $p\in G$ and where we consider functions from $G$ to $A$ as sitting inside $M(B\ot A)$ in the obvious way. Then $\Gamma$ is a left coaction making $A$ into a left $B$-comodule bi-algebra. 
\hfill$\square$
\einspr

Just remark that $\Gamma(A)(B\ot 1)\subseteq B\ot A$ because we have 
$$\Gamma(a)(\delta_p\ot 1)=\delta_p \ot \alpha_{p^{-1}}(a)$$
when $\delta_p$ is the function on $G$ that is $1$ in $p$ and $0$ everywhere else.
\snl
Because the algebra $B$ is abelian, we get a matched pair when we take the trivial right action of $A$ on $B$.

\iinspr{4.11} Proposition \rm
We assume the situation as in the previous proposition, but we now also assume that $(A,\Delta)$ has integrals. Then there is a homomorphism $\nu$ from $G$ to the non-zero complex numbers such that
$$\varphi(\alpha_p(a))=\nu(p)\varphi(a) \qquad\quad\text{and}\qquad\quad \psi(\alpha_p(a))=\nu(p)\varphi(a)$$
for all $p\in G$ and $a\in A$. In particular, the multiplier $y$, as obtained in Proposition 1.4 is given by the function $p\mapsto \nu(p^{-1})$.
\hfill$\square$
\einspr

Again, this is standard. Indeed, because $\alpha_p$ is assumed to be a coalgebra map, we must have that $\varphi\circ\alpha_p$ is again a left integral. By uniqueness of left integrals, there exists a complex number $\nu(p)$ so that $\varphi\circ\alpha_p=\nu(p)\varphi$. Because $\alpha$ is an action of $G$, we have that $\nu$ is a homomorphism. Finally, by the uniqueness of the antipode, we need to have that $\alpha_p\circ S=S\circ \alpha_p$ and as $S$ converts the left integral to a right integral, we also get $\psi\circ\alpha_p=\nu(p)\psi$. Compare this result with the one in Proposition 1.5.
\nl
Finally, we give an example where $y\neq 1$.

\iinspr{4.12} Example \rm Let $(A,\Delta)$ be an algebraic quantum group. Let $G$ be the additive group $\Bbb Z$ of integers. It acts on $(A,\Delta)$ by $\alpha_n(a)=S^{2n}(a)$. The scaling constant $\nu$ is a complex number such that $\varphi\circ S^2=\nu\varphi$. So the homomorphism in the previous proposition will be the map $n\mapsto \nu^n$. In particular, if the scaling constant differs from $1$, we get an example of a matched pair $(A,B)$ with a trivial right action of $A$ on $B$ and such that the multiplier $y$ is not equal to $1$.
\hfill$\square$\einspr
 
By a procedure as in Example 4.7, it is possible to construct cases where not only $y$ is non-trivial but also all the modular data are non-trivial. 
\nl
\it The special case where the left coaction of $B$ on $A$ is trivial \rm
\nl
Before we attack this case, first look again at the situation of the previous item. We assume that we have algebraic quantum groups so that we can look at the dual pair $(C,D)$ of the original pair $(A,B)$ as in Section 3. When the right action of $A$ on $B$ is trivial, then the right coaction of $C$ on $D$ is trivial. Indeed, we have (see Theorem 3.3)
$$\langle b\ot a,\Gamma(d)\rangle = \langle b\tl a, d\rangle =\varepsilon_A(a) \langle b,d\rangle $$
for all $a\in A$ and $b\in B$ so that $\Gamma(d)=d\ot 1$ for all $d\in D$. From the fact that $A$ is a left $B$-comodule bi-algebra here, it follows that the left action of $D$ on $C$ makes $C$ into a left $D$-module bi-algebra, so that
$$\Delta(d\tr c)=\sum_{(c)(d)}(d_{(1)}\tr c_{(1)}) \ot (d_{(2)}\tr c_{(2)})$$
for all $c\in C$ and $d\in D$ (to be interpreted with the necessary coverings). Finally, because  $T=T^{\text{op}}$
on $A\ot B$, we have that $R=R^{\text{op}}$ on $D\ot C$. This means that 
$$\sum_{(d)} (d_{(1)}\tr c) \ot d_{(2)} =\sum_{(d)} (d_{(2)}\tr c) \ot d_{(1)}$$
for all $c,d$ (see formulas (3.3) and (3.4) in the previous section). This equation will be satisfied if either the action of $D$ on $C$ is trivial (a case that corresponds to the coaction $\Gamma$ of $B$ on $A$ being trivial) or if $D$ is cocommutative (a case that is fulfilled if $B$ is abelian). 
\snl
This indicates what will happen when we take the case of multiplier Hopf algebras $A$ and $B$ where we now assume that the coaction $\Gamma$ of $B$ on $A$ is trivial. We get the dual version of the case of a trivial action.

\iinspr{4.13} Proposition \rm
Let $(A,B)$ be a matched pair and assume that the coaction $\Gamma$ of $B$ on $A$ is trivial, i.e.\ that $\Gamma(a)=1\ot a$ for all $a$ in $A$. Then we must have that $B$ is a right $A$-module {\it bi}-algebra and $R=R^{\text{op}}$, i.e.
$$\sum_{(a)}a_{(1)} \ot (b\tl a_{(2)})
= \sum_{(a)}a_{(2)}\ot (b\tl a_{(1)})$$
for all $a,b$. Conversely, if we have a right action of $A$ on $B$ making $B$ a right $A$-module {\it bi}-algebra and satisfying the above equality, then $(A,B)$ is a matched pair with the trivial coaction $B$ on $A$.
\hfill$\square$\einspr

Again, the proof is easy (and the result is stated already in [De-VD-W], see again the remark following Theorem 2.16 in that paper). Remark that the result follows by duality in the case of algebraic quantum groups, as explained above, but that in general, one still has to give an argument.
\snl
We will not look at more concrete examples in this case. Similarly as in Example 4.7, it is possible to construct cases where neither the action of $A$ is also trivial, neither $A$ has to be cocommutative.
\nl
\it The cases of a module and a comodule bi-algebra \rm
\nl
We have seen in one of the previous items that in the case of a trivial action of $A$ on $B$, it follows that $\Gamma$ will satisfy $\Gamma(aa')=\Gamma(a)\Gamma(a')$ for all $a,a'\in A$. This is of course also true for a trivial coaction - a case that we considered in a next item.
\snl
Now we will see what we can say if we (only) assume that $\Gamma$ is a homomorphism, i.e.\ when $A$ is assumed to be a $B$-comodule {\it bi}-algebra. We will give examples to show that this case is more general than any of the two previous cases where this happens.
\snl
Before we start with this investigation, we look again at the general case and we define the {\it fixed point subalgebra} in $M(B)$ of the right action of $A$ on $B$.
\snl
Recall that the action of $A$ on $B$ can be extended to the multiplier algebra $M(B)$ so that 
$$\align \sum_{(a)}(m\tl a_{(1)})(b\tl a_{(2)})&=(mb)\tl a \\
\sum_{(a)}(b\tl a_{(1)})(m\tl a_{(2)})&=(bm)\tl a
\endalign$$
for all $a\in A$, $b\in B$ and $m\in M(B)$. See e.g.\ Proposition 4.7 in [Dr-VD-Z]. Then, the following definition makes sense (see also Definition 4.10 in [Dr-VD-Z]).

\iinspr{4.14} Definition \rm
The fixed point algebra for the right action of $A$ on $B$ is defined as the set of elements $m\in M(B)$ satisfying $m\tl a=\varepsilon(a)m$ for all $a\in A$.
\hfill$\square$\einspr

It is not difficult to show that elements $m\in M(B)$ satisfying $m\tl a=\varepsilon(a)m$ for all $a\in A$ form a subalgebra in $M(B)$. One can also show that for these elements, we have that $ma=am$ in the multiplier algebra of the smash product $AB$ for all $a\in A$ and that also this property characterizes elements in the fixed point algebra. Finally remark that we define the fixed point algebra as a subalgebra of $M(B)$ and {\it not} as a subalgebra of $B$ itself (for obvious reasons).  
\snl 
Next we consider the 'left leg' of $\Gamma(A)$. Because it is assumed that both $\Gamma(a)(b\ot 1)$ and $(b\ot 1)\Gamma(a)$ belong to $B\ot A$ for all $a\in A$ and $b\in B$, we can define an element $(\iota \ot \omega)\Gamma(a)$ in $M(B)$ for all $a\in A$ and any linear functional $\omega$ on $A$. We can look at the space spanned by such elements in $M(B)$ for various choices of the space of linear functionals. For our purpose, we will consider the full dual space $A'$.
\snl
We can then show the following.

\iinspr{4.15} Proposition \rm
As before, let $(A,B)$ be a matched pair of multiplier Hopf algebras. Assume moreover that the coaction $\Gamma$ satisfies $\Gamma(aa')=\Gamma(a)\Gamma(a')$ for all $a,a'\in A$ (so that $A$ is $B$-comodule {\it bi}-algebra). Then the left leg of $\Gamma(A)$ is a subset of the fixed point algebra in $M(B)$ of the action of $A$ on $B$.

\snl\bf Proof\rm:
Take $a,a'\in A$ and rewrite the formula $\Delta_\#(aa')=\Delta_\#(a)\Delta_\#(a')$, using already that $\Gamma$ is multiplicative. This gives
$$\sum_{(a)(a')} (a_{(1)}a'_{(1)}\ot 1)\Gamma(a_{(2)})\Gamma(a'_{(2)})=
\sum_{(a)(a')} (a_{(1)}\ot 1)\Gamma(a_{(2)})(a'_{(1)}\ot 1)\Gamma(a'_{(2)}).$$
Replace in this formula $a$ by $a_{(2)}$ and multiply from the left in the first factor with $S(a_{(1)})$. Then the equation becomes
$$\sum_{(a')} (a'_{(1)}\ot 1)\Gamma(a)\Gamma(a'_{(2)})=
\sum_{(a')} \Gamma(a)(a'_{(1)}\ot 1)\Gamma(a'_{(2)}).$$
Now multiply from the right with $\Gamma(a'')$ where $a''\in A$ and replace $\sum_{(a')}a'_{(1)}\ot a'_{(2)}a''$ by $p\ot q$ where $p,q\in A$. Then we find
$$(p\ot 1)\Gamma(a)\Gamma(q)=\Gamma(a)(p\ot 1)\Gamma(q)$$
for all $a,p,q\in A$. Multiply from the right in the first factor with an element of $B$ and use that the cotwist map $T:q\ot b\mapsto \Gamma(q)(b\ot 1)$ is surjective. This shows that we can cancel the factor $\Gamma(q)$ in the previous formula, thus arriving at 
$$(p\ot 1)\Gamma(a)=\Gamma(a)(p\ot 1)$$
for all $a,p\in A$. This proves the result.
\hfill$\square$\einspr

We see that conversely, if $(p\ot 1)\Gamma(a)=\Gamma(a)(p\ot 1)$ holds for all $a,p\in A$, it will follow that $\Delta_\#$ is an algebra map on $A$ if this is the case also for $\Gamma$.
\snl
It seems that in this case, not much more can be said about the other conditions of a matched pair. Therefore, it is more interesting to consider now an example where the coaction $\Gamma$ is an algebra map, but where neither the coaction, nor the action  are trivial. We construct such an example with matched pairs of groups in the next item.
\snl
Before we do that, let us first say something about the case where it is assumed that $B$ is a right module-{\it bi}-algebra. As we mentioned already in Section 2 (see Remark 2.11), we will have that the integrals are relatively invariant under the action. This simplifies slightly the formulas for the modular automorphisms $\sigma_\#$ and $\sigma'_\#$ on $A$.
\nl
\it Back to matched pairs of groups \rm
\nl
Let us look at Example 6.2.16 in [M]. It is easy to generalize.  We take advantage of this to adapt the formulation to get it in accordance with our notations.
\snl
The starting point is any nilpotent ring $R$. So for every element $r\in R$ there exists an element $n\in \Bbb N$ satisfying $r^n=0$. Of course, we think here of the ring of $n\times n$ upper triangular matrices (say over the real numbers) with $0$ on the diagonal. Denote by $\widetilde R$ the ring obtained from $R$ by adding an identity and consider the subset $H$ of elements of the form $1 + r$ where $r\in R$. It is not difficult to show that $H$ is a group for the multiplication inherited from $\widetilde R$.
\snl
Next denote $K=H^{\text{op}}$, the group that is obtained by taking again $H$, but now with the opposite product. We will use $\gamma$ and $\gamma'$ for the identity maps, from $H$ to $K$ and from $K$ to $H$ respectively. These maps are each others inverses. And they are anti-isomorphisms.
\snl
Finally define
$$\align \theta\,:&H\to K \qquad \text{by} \qquad \theta(h)=\gamma(h)^{-1} \\ 
\theta':&K\to H \qquad \text{by} \qquad \theta'(k)=\gamma'(k)^{-1}.
\endalign$$
These maps are isomorphisms.
\snl
Then we have the following result.

\iinspr{4.16} Proposition \rm Let $H$, $K$, $\theta$ and $\theta'$ be as above. Define 
$$\align h\tr k &= 1+\theta(h)(k-1) \\
h\tl k &= 1+(h-1)\theta'(k)\endalign$$
for $h\in H$ and $k\in K$. Then $\tr$ left action of $H$ on $K$ and $\tl$ is a right action of $K$ on $H$, making $(H,K)$ a matched pair of groups (as in the first item of this Section).

\snl\bf Proof\rm:
The proof is straightforward. One can e.g.\ verify that 
$$\theta(h)k=(h\tr k)\theta(h\tl k)$$
for all $h\in H$ and $k\in K$ and if this is applied with $kk'$, $k$ and $k'$, we find
$$\align
\theta(h)kk'&=(h\tr (kk'))\theta(h\tl (kk'))\\
\theta(h)kk'&=(h\tr k)\theta(h\tl k)k' \\
            &=(h\tr k)((h\tl k)\tr k')\theta((h\tl k)\tl k')
\endalign$$
and we see that 
$$
(h\tr (kk'))=(h\tr k)((h\tl k)\tr k')
$$
for all $h\in H$ and $k,k'\in K$.
\hfill$\square$\einspr

This matched pair gives rise to a bicrossproduct as explained in Example 4.1.
\snl
Let us now specialize again and take for $R$ the ring of $3\times 3$ upper triangular matrices (say over $\Bbb R$) with $0$ on the diagonal.
\snl
We arrive at the following example.

\iinspr{4.17} Example \rm
i) The group $H$ is isomorphic with the semi-direct product $\Bbb R^2 \times_\alpha \Bbb R$ of the additive group $\Bbb R^2$ with the (left) action $\alpha$ of $\Bbb R$ on $\Bbb R^2$ given by 
$$\alpha_{p}(q,r)=(q,r+pq)$$
where $p,q,r\in \Bbb R$. The isomorphism is given by 
$$(q,r;p)\mapsto \left( \matrix 1 & p & r \\ 0 & 1 & q \\ 0 & 0 & 1 \endmatrix \right).$$
Similarly, the group $K$ is isomorphic with the semi-direct product $\Bbb R _\beta\times \Bbb R^2$ of the additive group $\Bbb R^2$ with the  action $\beta$ of $\Bbb R$ on $\Bbb R^2$ given by the same formula as the left action $\alpha$, but now considered as right action.  The isomorphism is now given by 
$$(p;q,r)\mapsto \left( \matrix 1 & 0 & 0 \\ p & 1 & 0 \\ r & q & 1 \endmatrix \right).$$
\snl
ii) The left action $\tr$ of $H$ on $K$ is now given by
$$(q,r;p)\tr (p';q',r')=(p';q',r'-qp')$$
while the right action $\tl$ of $K$ on $H$ is given by
$$(q,r;p)\tl (p';q',r')=(q, r-q'p;p).$$
\hfill$\square$\einspr

This example can again be generalized to more general types of semi-direct products of groups, using the same pattern. We leave this as an exercise for the reader.
\snl
In the following result, we see what is peculiar about this example (which is the main reason for including it here).

\iinspr{4.18} Proposition \rm In the previous example, we have 
$$\align h\tl (kk') &= (h\tl k)(h\tl k') \\
(hh')\tr k &= (h\tr k)(h'\tr k)\endalign$$
for all $h,h'\in H$ and $k,k'\in K$. In particular, this matched pair of groups gives rise to a matched pair of multiplier Hopf algebras $(A,B)$ (as in Theorem 2.14 and Theorem 3.1 in [De-VD-W]) so that the right action of $A$ on $B$ makes $B$ into a $A$-module {\it bi}-algebra and the left coaction of $B$ on $A$ makes $A$ into a $B$-comodule {\it bi}-algebra. Still, neither the action, nor the coaction are trivial.

\snl\bf Proof\rm: The proof is again straightforward. If we consider e.g.\ the left action $\tr$ of $H$ on $K$, we see that it is only the $\Bbb R^2$-component of $H$ that acts non-trivially on $K$ and one verifies that it is actually an action by automorphisms of $K$.
 \hfill$\square$\einspr
 
 It is clear that many more interesting examples can be obtained. Because we are working with multiplier Hopf algebras, the setting is less restrictive. In fact, in [L-VD] we obtain some more examples involving the p-adic numbers. Such examples can not be constructed within the setting of Hopf algebras. 
\nl\nl

\bf 5. Conclusions and further research \rm
\nl
In the present paper, we have continued our work on bicrossproducts of multiplier Hopf algebras ([De-VD-W]). In stead of working with matched pairs of general multiplier Hopf algebras, we  now took multiplier Hopf algebras with integrals (algebraic quantum groups). We showed already in [De-VD-W] that in that case, also the bicrossproduct has integrals. 
\snl
In [De-VD-W] we have given a formula for the right integral on the bicrossproduct in terms a right integrals on the components. In this paper we were able to give also the formula for the left integral, as well as the modular data: the modular element, the modular automorphisms, the scaling constant. This has been done in Section 2, using the preliminary results obtained in Section 1.
\snl
In Section 3 of this paper, we gave a satisfactory treatment of duality for the bicrossproducts. And in Section 4, we illustrated many of results by the use of examples.
\snl
Further research on this topic could focus on finding more examples, in particular examples that are typical for multiplier Hopf algebras and that are not obtainable as merely adaptations of known examples with Hopf algebras.
\snl
We refer to the work in progress on bicrossproducts of totally disconnected groups [L-VD]. In that paper, we consider also so-called {\it non-genuine} matched pairs of (locally compact) groups, mostly totally disconnected groups. It seems to be an interesting question if there exists a more general notion of matched pairs ('non-genuine matched pairs of multiplier Hopf algebras') so that the example of non-genuine matched pairs of groups (as studied in [L-VD]) fits into this more general theory.
\snl
Finally, it would be nice if some applications of our theory could be developed.
\nl\nl

\bf References \rm
\nl
{\bf [B-*]} M.\ Beattie, S.\ D\u{a}sc\u{a}lescu, L. Gr\"unenfelder \& C.\ N\u{a}st\u{a}sescu: {\it Finiteness conditions, co-Frobenius Hopf algebras and quantum groups}. J.\ Algebra {\bf 200} (1998), 312-333.
\snl
{\bf [DC-VD]} K.\ De Commer \& A.\ Van Daele: {\it Multiplier Hopf algebras embedded in locally compact quantum groups}.  Rocky Mountain Journal of Mathematics {\it 40} (2010), 1149-1182.
\snl
{\bf [De1]} L.\ Delvaux: {Semi-direct products of multiplier Hopf algebras: smash products}. Commun.\ Algebra {\bf 30} (2002), 5961-5977. 
\snl
{\bf [De2]} L.\ Delvaux: {Semi-direct products of multiplier Hopf algebras: smash coproducts}. Commun.\ Algebra {\bf 30} (2002), 5979-5997. 
\snl
{\bf [De-VD1]} L.\ Delvaux \& A.\ Van Daele: {\it The Drinfel'd double of multiplier Hopf algebras}. J.\ Algebra {\bf 272} (2004), 273-291.
\snl
{\bf [De-VD2]} L.\ Delvaux \& A.\ Van Daele: {\it Algebraic quantum hypergroups}. Adv.\ Math.\ {\bf 226} (2011), 1134-1167.
\snl
{\bf [De-VD-W]} L.\ Delvaux, A.\ Van Daele \& S.\ Wang: {\it Bicrossproducts of multiplier Hopf algebras}.  J.\ Algebra {\bf 343} (2011), 11-36. See Arxiv math.RA/0903.2974 for an expanded version.
\snl
{\bf [Dr-VD]} B.\ Drabant \& A.\ Van Daele: {\it Pairing and quantum double of multiplier Hopf algebras}. Algebras and Representation Theory {\bf 4} (2001), 109-132.
\snl
{\bf [Dr-VD-Z]} B.\ Drabant, A.\ Van Daele \& Y.\ Zhang: {\it Actions of multiplier Hopf algebras}. Commun.\ Algebra {\bf 22} (1999), 4117-4127.
\snl
{\bf [L-VD]} M.B.\ Landstad \& A.\ Van Daele: {\it Bicrossproducts for totally disconnected groups}. Preprint University of Trondheim and University of Leuven (2012). In preparation.
\snl
{\bf[M]} S.\ Majid: {\it Foundations of quantum group theory}. Cambridge University Press (1995).
\snl
{\bf [VD1]} A.\ Van Daele: {\it Multiplier Hopf algebras}. Trans.\ Amer.\ Math.\ Soc.\ {\bf 342} (1994), 917-932.
\snl
{\bf [VD2]} A.\ Van Daele: {\it An algebraic framework for group duality}. Adv.\ Math.\ {\bf 140} (1998), 323-366.
\snl
{\bf [VD3]} A.\ Van Daele: {\it Discrete quantum groups}. J.\ Algebra {\bf 180} (1996), 799-850.  
\snl
{\bf [VD4]} A.\ Van Daele: {\it Tools for working with multiplier Hopf algebras}. Arabian Journal for Science and Engineering {\bf 33} 2C (2008), 505--527. See also Arxiv math.RA/0806.2089. 
\snl
{\bf [VD-W]} A.\ Van Daele \& S.\ Wang: {\it A class of multiplier Hopf algebras}. Algebras and Representation Theory {\bf 10} (2007), 441-461.
\snl
{\bf [VD-Z]} A.\ Van Daele \& Y.\ Zhang: {\it A survey on multiplier Hopf algebras} In 'Hopf algebras and Quantum Groups', eds. S.\ Caenepeel \& F.\ Van Oyestayen, Dekker, New York (1998), pp. 259--309.
\snl

\end

{\bf [B-S]} S.\ Baaj and G.\ Skandalis: {\it Unitaires multiplicatifs et dualit\'e pour les produits crois\'es de C$^*$-alg\`ebres}. Ann.\ Scient.\ Ec.\ Norm.\ Sup.\ 4$^{\text{e}}$ s\'erie {\bf 26} (1993), 425--488. 
\snl
{\bf [B-*]} M.\ Beattie, S.\ D\v{a}sc\v{a}lescu, L.\ Gr\"unenfelder and C.\ N\v{a}st\v{a}sescu: {\it Finiteness conditions, co-Frobenius Hopf algebras and quantum groups}. J.\ Algebra {\bf 200} (1998), 312-333.
\snl
{\bf [B-M]} T.\ Brzezinski and S.\ Majid: {\it Quantum Geometry of Algebra Factorizations and Coalgebra Bundles}. Commun.\ Math.\ Phys.\ {\bf 213} (2000), 491-521.
\snl
{\bf [De1]} L.\ Delvaux: {\it Semi-direct products of multiplier Hopf algebras: Smash products}. Comm. Algebra {\bf 30} (2002), 5961--5977.
\snl
{\bf [De2]} L.\ Delvaux: {\it Semi-direct products of multiplier Hopf algebras: Smash coproducts}. Comm. Algebra {\bf 30} (2002), 5979-5997.
\snl
{\bf [De3]} L.\ Delvaux: {\it Twisted tensor product of multiplier Hopf ($^\ast$-)algebras}. J.\ Algebra {\bf 269} (2003), 285-316.
\snl
{\bf [De4]} L.\ Delvaux: {\it Twisted tensor coproduct of multiplier Hopf ($^\ast$-)algebras}. J.\ Algebra {\bf 274} (2004), 751--771.
\snl
{\bf [Dr-VD]} B.\ Drabant \& A. Van Daele: {\it Pairing and Quantum double of multiplier Hopf algebras}.  Algebras and Representation Theory 4 (2001), 109-132.
\snl
{\bf [Dr-VD-Z]} B.\ Drabant, A.\ Van Daele \& Y.\ Zhang: {\it Actions of multiplier Hopf algebras}. Comm.\ Algebra {\bf 27} (1999), 4117-4127.
\snl
{\bf [L-VD]} M.B.\ Landstad \& A.\ Van Daele: {\it Groups with compact open subgroups and multiplier Hopf $^*$-algebras}. Expo.\ Math.\ 26 (2008), 197--217.
\snl
{\bf [M1]} S.\ Majid: {\it Physics for algebraists: Non-commutative and non-cocommutative Hopf algebras by a bicrossproduct construction}. J.\ Algebra {\bf 130} (1990), 17-64.
\snl
{\bf [M2]} S.\ Majid:  {\it Hopf-von Neumann algebra bicrossproducts, Kac algebra bicrossproducts and the classical Yang-Baxter equations}. J.\ Funct.\ Analysis {\bf 95} (1991), 291-319.
\snl
{\bf [M3]} S.\ Majid: {\it Foundations of quantum group theory}. Cambridge University Press, 1995.
\snl
{\bf [V1]} S.\ Vaes: {\it Locally compact quantum groups}. Ph.D. thesis K.U.Leuven (2001).
\snl
{\bf [V2]} S.\ Vaes: {\it Examples of locally compact quantum groups through the bicrossed product construction}. Proceedings of the XIIIth Int.\ Conf.\ Math.\ Phys.\ London (2000). Editors A.\ Grigoryan, A.\ Fokas, T.\ Kibble and B.\ Zegarlinski, International press of Boston, Somerville MA (2001), pp.341--348.
\snl
{\bf [V-V1]} S.\ Vaes and L.\ Vainerman: {\it On low dimensional locally compact quantum groups}. Locally Compact Quantum Groups and Groupoids. Proceedings of the Meeting of Theoretical Physicists and Mathematicians, Strasbourg (2002). Ed.\ L.\ Vainerman, IRMA Lectures on Mathematics and Mathematical Physics, Walter de Gruyter, Berlin, New York (2003), pp.127-187.
\snl
{\bf [V-V2]} S.\ Vaes and L.\ Vainerman: {\it Extensions of locally compact quantum groups and the bicrossed product construction}. Adv.\ in Math.\ {\bf 175} (2003), 1--101.  
\snl
{\bf [VD3]} A.\ Van Daele: {\it Multiplier Hopf $^*$-algebras with positive integrals: A laboratory for locally compact quantum groups}. Locally Compact Quantum Groups and Groupoids. Proceedings of the Meeting of Theoretical Physicists and Mathematicians, Strasbourg (2002). Ed.\ L.\ Vainerman, IRMA Lectures on Mathematics and Mathematical Physics 2, Walter de Gruyter, Berlin, New York (2003), pp.229-247.
\snl
{\bf [VD4]} A. Van Daele: {\it Multiplier Hopf algebras with integrals}. Talk at the AMS Sectional Meeting, March 18-19, 2005 (Bowling Green, Kentucky, USA).
\snl

{\bf [VD-VK]} A.\ Van Daele \& S.\ Van Keer: {\it The Yang-Baxter and the Pentagon equation}. Comp.\ Math.\ {\bf 91} (1994), 201--221.
\snl 
{\bf [VD-W]} A.\ Van Daele \& S.\ Wang: {\it A class of multiplier Hopf algebras}. Algebras and Representation Theory {\bf 10} (2007), 441--461.
\snl

{\bf [VD-Z2]} A.\ Van Daele and Y.\ Zhang: {\it Galois Theory for multiplier Hopf algebras with integrals}. Algebras and Representation Theory {\bf 2} (1999), 83-106.
\end

-------------------------------------------------------------------------------------------------------------------------

\snl{\bfDC-VD]} K.\ De Commer \& A.\ Van Daele: {\it Multiplier Hopf algebras embedded in C$^*$-algebraic quantum groups}. Preprint K.U.\ Leuven (2006).
\snl{\bfDe-VD]} L.\ Delvaux \& A. Van Daele: {\it Algebraic quantum hypergroups}. Preprint University of Hasselt and K.U.\ Leuven (2006). math.RA/0606466
\snl{\bfF]} G.B.\ Folland: {\it A course in abstract harmonic analysis}. Studies in advanced mathematics. CRC Press, Boca Raton, London (1995).
\snl{\bfK-V1]} J.\ Kustermans \& S.\ Vaes: {\it A simple definition for locally compact quantum groups}. C.R.\ Acad.\ Sci.\ Paris S\'er I 328 (1999), 871--876.
\snl{\bfK-V2]} J.\ Kustermans \& S.\ Vaes: {\it Locally compact quantum groups}. Ann.\ Sci.\ \'Ecole Norm.\ Sup. 33 (2000), 837--934.
\snl{\bfK-VD]} J.\ Kustermans \& A.\ Van Daele: {\it C$^*$-algebraic quantum groups arising from algebraic quantum groups}. Int.\ J.\ Math.\ 8 (1997), 1067--1139.
\snl{\bfVD1]} A.\ Van Daele: {\it Multiplier Hopf algebras}. Trans.\ Am.\ Math.\ Soc.\ 342 (1994), 917--932.
\snl{\bfVD2]} A.\ Van Daele:  {\it An algebraic framework for group duality}.  Adv.\ Math.\ 140 (1998), 323--366.
\snl{\bfVD3]} A.\ Van Daele: {\it Locally compact quantum groups. A von Neumann algebra approach}. Preprint K.U.\ Leuven (2006), math.OA/0602212
\snl{\bfVD-W]} A.\ Van Daele \& Shuanhong Wang: {\it The Larson-Sweedler theorem for multiplier Hopf algebras}. J.\ Alg.\ 296 (2006), 75--95.
\snl{\bfVD-Z]} A.\ Van Daele \& Y.\ Zhang: {\it A survey on multiplier Hopf algebras} In 'Hopf algebras and Quantum Groups', eds. S.\ Caenepeel \& F.\ Van Oyestayen, Dekker, New York (1998), pp. 259--309.

\end

\snl{\bfB-B-T]} M.\ Beattie, D.\ Bulacu \& B.\ Torrecillas: {\it Radford's $S^4$ formula for co-Frobenius Hopf algebras}. J.\ Algebra (2006), to appear. 
\snl{\bfD]} L.\ Delvaux: {\it The size of the intrinsic group of a multiplier Hopf algebra}. Commun.\ Alg.\ 31 (2003), 1499--1514.
\snl{\bfD-VD1]} L.\ Delvaux \& A.\ Van Daele: {\it Traces on (group-cograded) multiplier Hopf algebras}. Preprint University of Hasselt and K.U.\ Leuven (2005).
\snl{\bfD-VD3]} L.\ Delvaux \& A.\ Van Daele: {\it The Drinfel'd double of multiplier Hopf algebras}. J.\ Alg.\ 272 (2004), 273--291.
\snl{\bfD-VD-W1]} L.\ Delvaux, A. Van Daele \& Shuanhong Wang: {\it Bicrossed product of multiplier Hopf algebras}. Preprint University of Hasselt, K.U.\ Leuven and Southeast University Nanjing (2005).
\snl{\bfD-VD-W2]} L.\ Delvaux, A. Van Daele \& Shuanhong Wang. {\it Quasi-triangular and ribbon multiplier Hopf algebras} Preprint University of Hasselt, K.U.\ Leuven and Southeast University Nanjing. In preparation.
\snl{\bfDr-VD]} B.\ Drabant \& A. Van Daele: {\it Pairing and Quantum double of multiplier Hopf algebras}.  Algebras and Representation Theory 4 (2001), 109-132.
\snl{\bfDr-VD-Z]} B.\ Drabant, A.\ Van Daele \& Y.\ Zhang: {\it Actions of multiplier Hopf algebras}. Commun.\ Alg.\ 27 (1991), 4117--4172.
\snl{\bfK]} J.\ Kustermans: {\it The analytic structure of algebraic quantum groups}. J.\ Algebra 259 (2003), 415--450.

\snl{\bfK-V3]} J.\ Kustermans \& S.\ Vaes: {\it Locally compact quantum groups in the von Neumann algebra setting}. Math.\ Scand.\ 92 (2003), 68--92. 

\snl{\bfK-VD2]} J.\ Kustermans \& A.\ Van Daele: {\it The Heisenberg commutation relations for an algebraic quantum group}. Preprint K.U.\ Leuven (2002), unpublished.
\snl{\bfL]} R.G.\ Larson: {\it Characters of Hopf algebras}. J.\ Algebra 17 (1971), 352--368.

\snl{\bfR]} D.\ Radford: {\it The order of the antipode of any finite-dimensional Hopf algebra is finite}. Amer.\ J.\ Math.\ 98 (1976), 333--355. 

\snl{\bfVD3]} A.\ Van Daele: {\it Multiplier Hopf algebras with integrals}. Talk at the AMS meeting in Bowling Green, Kentucky (2005).

\snl{\bfVD-W]} A.\ Van Daele \& Shuanhong Wang: {\it Modular multipliers of regular multiplier Hopf algebras of discrete type}. Preprint K.U.\ Leuven (2004), unpublished.

\end

M. Beattie, S. D$\check{a}$sc$\check{a}$lescu, L.
 Gr$\ddot{u}$nenfelder, C. N$\check{a}$st$\check{a}$sescu,
 Finiteness conditions, Co-Frobenius Hopf algebras
 and Quantum groups, J.\ Alg.\ {\bf 200} (1998), 312--333.
\snl

[De1] L. Delvaux, The size of the intrinsic group of a multiplier
 Hopf algebra. Comm. in Algbera, 31(3)(2003), 1499-1514.
\snl
[De2] L. Delvaux, Twisted tensor coproduct of multiplier Hopf
 algebras, J. Algbera, 2004, to appear
\snl
[De3] L. Delvaux, Pairing and Drinfel'd double of Ore-extensions,
 Comm. in Algebra 29(7)(2001), 3167-3177.
\snl

[De-VD] L. Delvaux and A. Van Daele, The Drinfel'd double versus
 the Heisenberg double for an algebraic quantum group,
 J. Pure and Appl. Algbera, 2004, to appear.
\snl

[D-VD] B. Drabant and A. Van Daele, Pairing and the quantum double
 of multiplier Hopf algebras, Algebras and Representation Theory
 4:(2001), 109-132.
\snl
[D-VD-Z] B. Drabant,  A. Van Daele, and Y. Zhang, Actions
 of multiplier Hopf algebras, Comm. Algebra, 27(9)(1999),
 4117-4172.
\snl
[KR] L. H. Kauffman and D. E. Radford, A necessary and sufficient
 condition for a finite-dimensional Drinfel'd double to be
 a ribbon Hopf algebra, J. Algebra, 159(1993), 98-114.
\snl

[Ma] S. Majid, Quasitriangular Hopf algebras and Yang-Baxter
 equations, Internat. J. Modern Phys. A 5(1990), 1-91.
\snl
[Pa] F. Panaite, Ribbon and charmed elements for quasitriangular
 Hopf algebras, Communications in algebra, 25(3)(1997), 973-977.
\snl
[R1] D. E. Radford, On the antipode of a quasitriangular Hopf
algebra, J. Algebra
 151(1992), 1-11.
\snl
[R2] D. E. Radford, The trace function and Hopf algebras, J.
Algebra 163(1994),
 583-622.
\snl
[RT] N. Yu. Reshetikhin and V. G. Turaev, Ribbon graphs and their
 invariants derived from quantum groups, Comm. Math. Phys.
 127(1990), 1-26
\snl
[VD-Z] A. Van Daele and Y. Zhang, Multiplier Hopf algebras of
 discrete type, J. Algebra, 214(1999), 400-417.

\end